\documentclass{amsart}
\usepackage{amssymb,amsmath,amscd}
\usepackage{epstopdf}
\usepackage{graphicx}

\newtheorem{thm}{Theorem}[section]

\newtheorem{corollary}[thm]{Corollary}
\newtheorem{prop}[thm]{Proposition}
\newtheorem{lemma}[thm]{Lemma}
\newtheorem{fact}[thm]{Fact}
\newtheorem{question}[thm]{Question}

\theoremstyle{definition}
\newtheorem{defn}[thm]{Definition}
\newtheorem{example}[thm]{Example}

\theoremstyle{remark}
\newtheorem{remark}[thm]{Remark}

\newcommand{\bt}{\begin{thm}}
\newcommand{\et}{\end{thm}}
\newcommand{\bp}{\begin{prop}}
\newcommand{\ep}{\end{prop}}
\newcommand{\bd}{\begin{defn}}
\newcommand{\ed}{\end{defn}}
\newcommand{\bl}{\begin{lemma}}
\newcommand{\el}{\end{lemma}}
\newcommand{\bfa}{\begin{fact}}
\newcommand{\efa}{\end{fact}}
\newcommand{\bc}{\begin{corollary}}
\newcommand{\ec}{\end{corollary}}
\newcommand{\bex}{\begin{example}}
\newcommand{\eex}{\end{example}}
\newcommand{\br}{\begin{remark}}
\newcommand{\er}{\end{remark}}
\newcommand{\ben}{\begin{enumerate}}
\newcommand{\een}{\end{enumerate}}

\newcommand{\sotto}[2]{#1_{#2}}

\newcommand{\rrr}{\rightarrow}
\newcommand{\ra}{\rightarrow}

\newcommand{\ds}{\displaystyle}

\newcommand{\ideal}[1]{\sotto {{\mathcal I}}{#1}}

\newcommand{\exact}[3]
{0 \rrr #1 \rrr #2
\rrr #3 \rrr 0}

\newcommand{\lexact}[3]
{0\rrr #1 \rrr #2
\rrr #3 }

\newcommand{\pso}{\mathbb{P}^3}

\newcommand{\PP}{\mathbb{P}}

\newcommand{\Z}{\mathbb{Z}}

\newcommand{\N}{\mathbb{N}}

\newcommand{\coo}{{\mathcal O}}

\newcommand{\cae}{{\mathcal E}}

\newcommand{\gon}{\mathrm{gon}}

\begin{document}

\title[Gonality of a general ACM curve]{Gonality of a general ACM curve in $\PP^3$.}
\author{R. Hartshorne and E. Schlesinger}

\address{Department of Mathematics, University  of
California, Berkeley, California 94720--3840}

\address{Dipartimento di Matematica, Politecnico di Milano, Piazza Leonardo da
Vinci 32, 20133 Milano, Italia}

\thanks{The first author was partially supported by Gnsaga - Programma Professori visitatori.
The second author was partially supported by
MIUR PRIN 2005:
  {\em Spazi di moduli e teoria di Lie}.
 }


\begin{abstract}
Let $C$ be an ACM (projectively normal) nonsingular curve in $\PP^3_{\mathbb{C}}$ not contained in a plane,
and suppose $C$ is general in its Hilbert scheme - this is irreducible once the postulation
is fixed.
Answering a question posed by Peskine, we show the gonality of $C$ is $d-l$, where $d$ is the degree
of the curve, and $l$ is the maximum order of a multisecant line of $C$.
Furthermore $l=4$ except for two series of cases, in which the postulation of $C$ forces every surface of
minimum degree containing $C$ to contain a line as well. We compute the value of $l$ in terms of the postulation of $C$
in these exceptional cases. We also show the Clifford index of $C$ is equal to $\gon(C)-2$.
\end{abstract}

\maketitle

\section{Introduction}
Let $C$ be a nonsingular projective curve over an algebraically closed field $\mathbb{K}$.
The {\em gonality}  of $C$, written $\gon (C)$, is the minimum
degree of a surjective morphism $C \ra \PP^1$, or equivalently the minimum positive integer
$k$ such that there exists a $g^1_k$ on $C$.

For curves of genus $g \geq 1$ the gonality varies between $2$, the value it takes on hyperelliptic curves,
and $\left[ \frac{g+3}{2}\right]$, which by Brill-Noether theory is the gonality of a general curve of genus $g$.
It may be regarded as the most fundamental invariant of the algebraic structure of $C$ after the genus,
providing a stratification of the moduli space of curves of genus $g$.

When a curve is embedded in some projective space, it is natural to wonder whether the gonality may be
related to extrinsic properties of the curve. A classical result in this direction, already known to Noether
- cf. \cite{Cil, HNoe}- is

\bt
A smooth curve $C \subset \PP^2$ of degree $d \geq 3$ has gonality $\gon (C)= d\!-\!1$, and
any morphism  $C \ra \PP^1$ of degree $d\!-\!1$ is obtained projecting $C$ from one
of its points.
\et
See \cite{mjm} for a proof and references. It is a simple exercise to prove the statement using Lazarsfeld's method \cite{LK3}
that associates a vector bundle on $\PP^2$ to a base point free pencil on $C$. It is this method that we will exploit
in the proof of our result.

One may ask a similar question for a curve  $C \subset \PP^3$. If $L$ is a line in $\PP^3$,
projection from $L$ induces a morphism $\pi_L: C \ra \PP^1$, whose degree is the
the degree of $C$ minus the number of points of intersection of $C$ and $L$. Thus the morphisms $\pi_L$ of minimal degree
are those corresponding to maximal order multisecant lines. We define
$$
l=l(C) = \mbox{Max} \{\deg (C \cap L): \mbox{$L$ a line in $\PP^3$} \}
$$

By analogy with the plane curves case one might wonder whether
\begin{equation}\label{question}
\gon(C)= \deg (C)-l(C)
\end{equation}
for a curve in $\PP^3$, in which case following the terminology of \cite{mjm} we say the gonality of $C$ is computed
by multisecants. Of course, this is usually not the case. For example, a general curve of genus $g$ has $\gon (C)=\left[ \frac{g+3}{2}\right]$
and can be embedded in $\PP^3$ as a nonspecial linearly normal curve of degree $g+3$. Since the Grassmannian of lines
in $\PP^3$ has dimension $4$, and the set of lines meeting $C$ is a codimension one subvariety, one expects
$l(C)$ to be $4$, and so
$$
\deg (C) - l(C) =g-1 >\left[ \frac{g+3}{2}\right]=\gon(C).
$$
See \cite[Examples 2.8 and 2.9]{mjm} for specific counterexamples.

On the other hand, if the embedding of $C$ in $\PP^3$ is very special, one may hope the gonality of $C$ is
computed by multisecants. In this vein Peskine raised the question:

\begin{question}
If $C$ is a smooth ACM curve in $\PP^3$, is its gonality computed by multisecants ?
\end{question}
Here ACM means arithmetically Cohen-Macaulay, and a curve in $\PP^3$ is ACM if and
only if the natural maps
$$
H^0 (\PP^3, \coo (n)) \ra H^0 (C, \coo_C (n))
$$
are surjective for every $n \geq 0$.

Some special cases have been treated in the literature. Early results about uniqueness of the linear series
$|\coo_C(1)|$ for complete intersections and other ACM curves are in \cite{CilL}.
Basili \cite{Bas}
has proven the gonality of a smooth complete intersection is indeed computed by multisecants; besides,
Ellia and Franco \cite{ef} showed that the maximum order $l$ of a multisecant to a general complete intersection
of type $(a,b)$ is $4$ if $a \geq b \geq 4$ as one expects. Lazarsfeld in \cite[4.12]{L} finds lower bounds for the
gonality of a complete intersection curve in $\PP^n$.

Results by Martens \cite{martens} and Ballico \cite{ballico} show that the gonality of a smooth curve
$C \subset \PP^3$ on a smooth quadric surface  is computed by multisecants. Hartshorne in \cite{mjm} shows that if
a smooth curve $C \subset \PP^3$ is ACM,  lies on a smooth cubic surface $X$, and is general in its linear system on $X$,
then  its gonality is computed by multisecants. Farkas \cite{farkas} has shown that  smooth ACM curves
$C \subset \PP^3$  lying on certain smooth quartic surfaces that do not contain rational or elliptic curves
have gonality computed by multisecants.

In this paper, we show that, with the exception of very few cases we cannot decide, the gonality of a
{\em general} ACM curve is indeed computed by multisecants.
We have to make sense of the expression general ACM curve. To obtain an irreducible parameter space for
ACM curves one needs to fix the {\em Hilbert function}, that is,
the sequence of integers $h^0(\coo_C(n))$. This is more conveniently expressed
by its second difference or $h$-vector:
$$
h_C (n)=
h^0(\coo_C (n))-2 h^0(\coo_C(n-1))+h^0(\coo_C(n-2)).
$$
which has the advantage of being finitely supported while still nonnegative. We will denote by
$A(h)$ the Hilbert scheme parametrizing ACM curves in
$\PP^3$ with $h$-vector $h$. By a theorem due to Ellingsrud \cite{E} (cf. \cite[p. 5, Corollaire 1.2 p. 134 and 1.7 p. 139]{MDP}),
the Hilbert scheme $A(h)$ is {\em smooth and irreducible}.
Thus by a general ACM curve we will mean a curve in a Zariski open nonempty subset of $A (h)$.
We believe it is reasonable to assume that $C$ is general in the statement of our theorem,
because it might happen that a special ACM curve had a low degree pencil unrelated to the line bundle $\coo_C(1)$.

\bt \label{st1} \
Assume $\mathbb{K}=\mathbb{C}$ is the field of complex numbers.
Let $C \subset \PP^3$ be a nonplanar smooth ACM curve. If $C$ is general in the Hilbert scheme $A(h_C)$, then
$$
\gon (C)= d - l
$$
where $d= \deg (C)$ and $l=l(C)$ is the maximum order of a multisecant line to $C$,
except perhaps if one of the following occurs (where $s$ denotes the least degree of a surface containing $C$):
\begin{itemize}
    \item $s=5$ and $(d,g)=(15,26),(16,30)$
    \item $s=6$ and $(d,g)=(21,50),(22,55),(23,60)$
    \item $s=7$ and $(d,g)=(28,85),(29,91)$
    \item $s=8$ and $(d,g)=(36,133)$
\end{itemize}
\et
\br
For curves $C$ contained in a quadric or a cubic surface, the statement follows from the references
cited above. So our contribution is for curves not lying on a cubic surface.
\er

We can also determine the integer $l(C)$ in terms of the $h$-vector of $C$. Most of the time
$l(C)=4$, with two families of exceptions. These exceptional cases arise because the $h$-vector
forces  surfaces of minimal degree containing $C$ to contain a line as well; this line
is then a multisecant of order higher than expected.

Denote by
\begin{itemize}
    \item $s$ the {\em least degree} of a surface containing $C$
    \item $t$ the integer
    $$\mbox{Min} \{n: h^0 (\ideal{C}(n))-h^0 (\coo_{\pso}(n-s))>0 \}$$
    \item $e$  the  {\em index of speciality} of $C$:
    $$e = \mbox{Max}\{n: h^1 \coo_C (n) >0\}$$
\end{itemize}

The value of $l(C)$ is given by:
\bt \label{st2}
Let $C \subset \PP^3_{\mathbb{C}}$ be a general smooth ACM curve with $s \geq 4$. Let $l=l(C)$
denote the maximum order of a multisecant line to $C$. Then:
\begin{itemize}
    \item $l=4$, unless
    \item the $h$-vector of $C$ satisfies $h(e+1)=3$ and $h(e+2)=2$, in which case
    $l=e+3$ and $C$ has a unique $(e+3)$-secant line,
    or
    \item $t>s+3$ and the $h$-vector of $C$ satisfies $h(t)=s-2$ and $h(t+1)=s\!-\!3$, but not $h(e+1)=3, h(e+2)=2$,
    in which case $l=t\!-\!s\!+\!1$ and $C$ has a unique $(t\!-\!s\!+\!1)$-secant line.
\end{itemize}
\et
Nollet \cite{nollet} has found a sharp bound for the maximal order $l=l(C)$
of a multisecant line in terms of the $h$-vector of $C$, valid for {\em any} irreducible ACM curve. If $C$ is not a complete intersection,
the bound is the largest integer $n$ for which
$$
h_C(n-1)-h_C(n) >1.
$$
He also shows there exist smooth curves in the Hilbert scheme $A(h)$ achieving the bound. As
this number is  $ \geq s_C$, we see that $l(C)$ and the gonality of $C$ {\em vary} in the family $A(h)$, provided
$s \geq 5$ and the gonality of the general curve is $d-4$ (in fact the argument of Theorem \ref{no5secants} shows that $l(C)$ varies
in the linear system $|C|$ on a smooth surface $X$ of degree $s=s_C \geq 5$ containing $C$).
On the other hand, in the special case $h(e+1)=3$ and $h(e+2)=2$,
then Nollet's bound is precisely $e+3$, so that $l(C)$ is constant in $A(h)$.

\vspace{.3cm}
Finally, in most cases we can prove that {\em every} pencil computing the gonality of $C$ arises
from a maximum order multisecant:
\bt \label{st3} \
Let $C \subset  \PP^3_{\mathbb{C}}$ be a general smooth ACM curve with $s \geq 4$.
Then
every pencil of minimal degree on $C$ arises from a maximal order multisecant line,
except
\begin{itemize}
    \item[a)] if one of the following cases we cannot decide occurs:
  \begin{itemize}
    \item  $s=4$ and $(d,g)=(10,11), (11,14),(12,17)$, or
    \item  $s=5$ and $(d,g)=(15,26),(16,30),(17,34),(18,38)$, or
    \item  $s=6$ and $(d,g)=(21,50),(22,55),(23,60),(24,65)$, or
    \item  $s=7$ and $(d,g)=(28,85),(29,91), (30,97)$, or
    \item  $s=8$ and $(d,g)=(36,133),(37,140)$; or if
  \end{itemize}
    \item[b)] $C$ is linearly equivalent to $C_0+bH$ on a smooth quartic surface,
    where $C_0$ is an elliptic quartic curve, $H$ is a plane section and $b \geq 2$.
 In this case the gonality is $d-4$, and $\coo_C(b)$ is a  $g^1_{d-4}$ on $C$ that does not arise from a $4$-secant.
\end{itemize}
In particular, $C$ has a finite number of pencils of minimal degree, and therefore its Clifford index  is
$$
\mbox{Cliff}(C)= \gon (C) -2= d-l(C)-2
$$
\et

We begin the paper illustrating the proof of the main theorem with two specific examples.
The paper is then structured according to the following outline of the proof. Since the conclusions of our result are semicontinuous
on the Hilbert scheme $A(h)$, it suffices to show the existence of a {\em single} curve $C$ for which the result holds.
Let $C$ be a smooth ACM curve in $\PP^3$ with given $h$-vector $h$,
not lying on any surface of degree $\leq 3$.
In section \ref{Cayley} we review the classical result that for every smooth space curve $D$ of degree $\geq 10$ there exists a line $L$ that is at
least a $4$-secant line of $D$. Thus $\gon (C) \leq d-4$.
Next, if $C$ is  general in $A(h)$, it is contained in
a {\it smooth} surface $X$ of degree $s$, where $s$ is the least degree of a surface containing $C$.
We prove in Corollary \ref{no5s} that, if $C$ is general in its linear system on $X$ and $L$ is an $l$-secant line of $C$ with $l \geq 5$,
then $L$ is contained in $X$. In fact, we prove a slightly more general result,
which gives explicit conditions for a space curve  not to have $5$-secant lines:
\bt[see Theorem \ref{no5secants} ]
Let $C \subset \PP^3_{\mathbb{K}}$ be a curve
contained in an irreducible surface $X$ of degree $s$.
Suppose $C$ is a Cartier divisor on $X$, and
$$
\begin{cases}
H^0 (\PP^3, \ideal{C} (s\!-\!2)\,)= 0  \\
H^1 (\PP^3, \ideal{C} (m)\,)= 0 \mbox{\ \ \ \ for $m=s\!-\!2, s\!-\!3, s\!-\!4$}
\end{cases}
$$
If $C$ is general in its linear system on $X$,
then $\deg (C \cap L) \leq 4$  for every line $L$ not contained in $X$, and $C$ has only finitely many $4$-secant
lines not contained in $X$.

In particular, if $X$ does not contain a line, then $C$ does not have an $l$-secant line for any $l \geq 5$.
\et

At this point to prove our main theorem we need to show that every pencil of minimal degree arises from a multisecant
line. The proof uses Lazarsfeld's technique \cite{LK3} that associates to a base point free pencil on $C$ a vector
bundle $\mathcal{E}$ on the surface $X$ as explained in section \ref{laz}.

In section \ref{liaison} we review enough liaison theory for ACM curves to be able to show that
the Lazarsfeld bundle $\mathcal{E}$ satisfies
$\Delta (\mathcal{E})= c_1 (\mathcal{E})^2 - 4 c_2 (\mathcal{E}) >0$
with the exceptions listed in the statement of Theorem \ref{st1}.
It follows by Bogomolov's theorem \cite{Bog} that, if $char. (\mathbb{K})=0$, then
$\mathcal{E}$ is Bogomolov unstable.
Thus it has a destabilizing
divisor $A \in Pic (X)$, whose degree $x=A.H$  satisfies stringent numerical restrictions in terms of the intersection
numbers of $A^2$, $A.C$ and $C^2$- see the proof of Theorem \ref{gon}.

To use effectively these restraints we need to control the Picard group of $X$.
Here we use the hypothesis the ground field is $\mathbb{C}$ to be able to apply the Noether-Lefschetz type Theorem of Lopez \cite[II.3.1]{lopez}
or the more recent work of  Brevik and Nollet \cite{bs} to conclude the following: if $C$ is general in $A(h)$ and $X$ is very general
among surfaces of minimal degree containing $C$, then
 $Pic (X)$ is freely generated by $H$ and the irreducible components of a curve $\Gamma$ that is general among
curves minimally linked to $C$.
$\Gamma$ is a general ACM curve, but  it may not be irreducible. Thus we are led to establish a structure theorem for general ACM curves.
Section \ref{gacm}
is devoted to the proof of this result. It generalizes Gruson-Peskine's theorem \cite{GP} according to which the
general ACM curve in $A(h)$ is smooth and irreducible if $h$ is of decreasing type ("has no gaps") :
\bt[see Theorem \ref{dthm} ]
Let $A(h)$ denote  the Hilbert
scheme parametrizing ACM curves in $\PP^3_{\mathbb{K}}$ with $h$-vector $h$.
If $\Gamma$ is general in $A(h)$, then
$$\Gamma =D_1 \cup D_2 \cup \ldots \cup D_r$$
   where $r+1$ is the number of Gruson-Peskine gaps of $h$, and the $D_i$ are distinct smooth irreducible ACM curves whose
   $h$-vectors are determined by the gap decomposition of $h$ as explained in
   section \ref{gacm}. Furthermore,
   for every $1 \leq i_1 < i_2 <  \cdots < i_h \leq r$, the curve
   $$
D_{i_1} \cup D_{i_2} \cup \ldots \cup D_{i_h}
   $$
   is still ACM.
\et

Thus we can write the destabilizing divisor as $A=aH + \sum a_i D_i$. In the proof of the main Theorem \ref{gon},
using the fact that the curves $D_i$ and their unions are ACM, together
with the numerical restraints on $x=A.H$ we show $-s\!-\!1 \leq x <0$
(see Section \ref{examples} for a specific example).
We then play this inequality against the bounds of Corollary \ref{cbounds}, which are essentially upper bounds for the genus
of an ACM curve lying on $X$ in terms of the degree of the curve and of degree of $X$. In fact, these bounds are a refinement
of the bounds for the genus of an ACM curve proven by Gruson and Peskine in \cite{GP} (see Remark \ref{remark}).
The end result is that there are only two possibilities for $A$:
either $-A=H$ (the plane section) or $-A=H-L$ for some line $L$ on $X$.

Corollary \ref{ms} shows that in case $A=-H$
the pencil arises from a multisecant line not contained in $X$, while in case $A=L-H$ the
pencil arises from $L$. This shows pencils of minimal degree on $C$ all arise from multisecant lines,
thus completing the proof of the theorem.

The second named author would like to thank Gian Pietro Pirola who explained to him Lazarsfeld's results on linear series
while working on \cite{PS}, and Cecilia Rizzi for several conversation.

\section{Two examples} \label{examples}
Before plunging into the general case,  we illustrate the proof in two specific examples.
\subsection{Example with $\Gamma$ reducible, but no line on $X$} \label{example1} \

Let $C$ be a general ACM curve with
$h$-vector $\{1,2,3,4,5,6,7,4,3\}$. Then $C$ has
$d=35$, $g=130$, $s=7$ and $e=6$. We claim the gonality of such a curve is $d-4=31$.

The curve $C$  can be
linked by two degree $7$ surfaces to a curve $\Gamma$ with
$$
h_\Gamma= \{1,2,3,4,2,2\}
$$

According to Theorem \ref{dthm}, the general such $\Gamma$ is the union of a twisted cubic curve $D_1=T$ and a curve $D_2=D$
of type $(6,5)$ on a smooth quadric.

Working over $\mathbb{C}$, we can assume by the Noether-Lefschetz type result of Lopez \cite{lopez} that there
exists such a pair $(C, \Gamma)$ on a smooth
surface $X$ of degree $7$ whose Picard group  $Pic(X)$ is freely generated by the classes of $H$, $T$ and $D$.

Throughout the paper we will make use of the following bilinear form on $Pic (X)$:
$$
\phi(D,E)= (D.H) \, (E.H) - s \, (D.E)=
\det
\begin{bmatrix}
  D.H & H^2 \\
  D.E & E.H
\end{bmatrix}
$$
where $s=H^2= \deg (X)$. This is essentially the positive definite product on $Pic (X)/{\Z H}$ induced by
the intersection product: by the algebraic Hodge index
theorem $ \phi(D,D) \geq 0$  for any divisor $D$ on $X$, and $\phi(D,D)=0$
if and only if $D$ is numerically (hence linearly) equivalent to a multiple of $H$.
The associated quadratic form is
$$
\phi (D,D)= d^2+s(s-4)d-2s(g-1)
$$
which in the case at hand gives
$$
\begin{cases}
  \phi(T,T) = 86 \\
  \phi(D,D) = 86 \\
  \phi(C,C) = 154
\end{cases}
$$

The class of $C$ in $Pic(X)$ is $7H-T-D$, hence
$$
154= \phi(C,C)= \phi(T+D,T+D)= 172+ 2 \phi(T,D)
$$
and we conclude $\phi(T,D)=-9$
(Proposition \ref{recursion}a allows to perform
this calculation in general: see formula (\ref{dij}) in the proof of  Theorem \ref{gon}).

Suppose $\mathcal{Z}$ is a complete base point free $g^1_k$ on $C$ with $k \leq 31$. Then
$C^2-4k=153-4k >0$, so the bundle associated to $\mathcal{Z}$ on $X$ is Bogomolov unstable
(see section \ref{laz}),
and comes with a destabilizing divisor $A$.

Write $A=aH+bT+cD$. Then
$$\phi(A,A)=\phi(bT+cD,bT+cT)=86b^2+86c^2-18bc$$
and
$$\phi(A,C)=\phi(bT+cD, -T-D)=-86b+9c+9b-86c=-77(b+c).$$
Therefore
$$
\phi(A,A+C) =77b(b-1)+77c(c-1)+9(b-c)^2.
$$

The fact that $A$ is the destabilizing divisor for the vector bundle associated to the pencil $\mathcal{Z}$
implies that its degree $x$ satisfies the following stringent
numerical restraints (see the proof of Theorem \ref{gon})
$$
\begin{cases}
    -35/2 < x < 0 \\
     x^2\geq \phi(A,A) \\
x^2+ 35 x +7k \geq 77b(b-1)+77c(c+1)+18(b-c)^2
\end{cases}
$$
From the last inequality we see $x^2+35x +217 \geq 0$
because $k \leq 31$. This together with $x<0$ forces $x \geq -8$.
Note that here we have used the fact that
$\phi(A, A+C) \geq 0$ to conclude $x \geq -8$. In the general case,
we will show that we still have $\phi(A, A+C) \geq 0$,
and that this allows to conclude $x \geq -s\!-\!1$ as in this example.

Unless $b=c=0$, the inequality $$x^2 \geq \phi(A,A)= 77b^2+77c^2+ 9(b-c)^2$$
gives $x^2 \geq 77$, which is absurd. Corollary \ref{s1bound} of Section \ref{qsection}
essentially shows that $\phi(A,A) > (s\!+\!1)^2$ also holds in general,
except for two series of cases which arise when $X$ contains a line.

We conclude $A=aH$, and then $-8 \leq x=7a \leq -1$ implies $A=-H$.
Then Corollary \ref{ms} shows that
the given pencil arises from a $(d-k)$-secant line not contained in $X$.
Finally by Corollary \ref{no5s}, if $C$ is general in its linear system on $X$ and $l \geq 5$,
then $C$ has no $l$-secant line that is not contained in $X$. Hence $k=d-4=31$.

\subsection{Example with $\Gamma$ reducible and a line on $X$} \

Let $C$ be a general ACM curve with
$h$-vector $\{1,2,3,4,5,3,2\}$.
Then $C$ has
$d=20$, $g=48$, $s=5$ and $e=4$. In this case $h_C(e+1)=3$ and $h_C(e+2)=2$, so our theorem
predicts that the gonality of $C$ is $13$ and that $C$ has a unique $g^1_{13}$, arising from
a $7$-secant line $L$. We explain the proof in this special case.

The curve $C$  can be
linked by two degree $5$ surfaces to a curve $\Gamma$ with $h$-vector $\{1,2,1,1\}$. Note that $d_\Gamma=5$, and  $g(\Gamma)=3$.
Since $C$ is general, $\Gamma$ can be taken general among ACM curves with $h$-vector $\{1,2,1,1\}$.
The general such curve is, according to \ref{dthm}, the union of a line $L$ and a plane quartic $P$
meeting at one point. As above we may assume the very general
quintic surface $X$ through $C$ is smooth and has $Pic(X) \cong \Z^3$ with
generators $H$, $L$ and $P$.

By the formula
$$
\phi (D,D)= d^2+s(s-4)d-2s(g-1)
$$
we can compute
$$
\begin{cases}
  \phi(L,L) = (s\!-\!1)^2= 16 \\
  \phi(P,P) = \phi(H-P, H-P)= 16 \\
  \phi(C,C) = \phi(L + P, L+ P)= 25+25-10 \times 2= 30
\end{cases}
$$

The class of $C$ in $Pic(X)$ is by construction $5H-L-P$, hence
$$
30= \phi(C,C)= \phi(L+P,L+P)= 32+ 2 \phi(L,P)
$$
and we conclude $\phi(L,P)=-1$.

It follows that $\phi(C,L)=\phi(-L-P,L)=-15$, hence $C.L=7$: the line $L$ is a $7$-secant
line of $C$, and therefore the gonality of $C$ is at most $13$.

We compute
$$C^2=\frac{d^2-\phi(C,C)}{5}= \frac{370}{5}=74$$
Suppose $\mathcal{Z}$ is a $g^1_k$ on $C$ with $k \leq 13$. Then
$C^2-4k=74-4k \geq 22 >0$, so the rank two bundle $\mathcal{E}$ associated to the given pencil is Bogomolov unstable,
and comes with a destabilizing divisor $A$.

Write $A=aH+bL+cP$. Then
$$\phi(A,A)=\phi(bL+cP,bL+cP)=16b^2+16c^2-2bc$$
and
$$\phi(A,C)=\phi(bL+c P, -L-P)=-16b+c+b-16c=-15(b+c).$$
Therefore
$$
\phi(A,A+C) =15b(b-1)+15c(c-1)+(b-c)^2
$$

Thus the inequalities for  $x= \deg(A)$ are in this case
$$
\begin{cases}
    -10 < x < 0 \\
     x^2\geq \phi(A,A) \\
x^2+20 x +5k \geq 15b(b-1)+15c(c-1)+(b-c)^2.
\end{cases}
$$
From the third inequality, keeping into account that $k \leq 13$, we see $x^2+20x +65 \geq 0$.
This together with $x<0$ forces $x \geq -4$.

If $-3 \leq x \leq -1$, then the second inequality gives
$$
9 \geq \phi(A,A)=16b^2+16c^2-2bc= 15b^2+15c^2+(b-c)^2
$$
hence $b=c=0$. But then $x= A.H=5a$, a contradiction.

Then $x=-4$, and since $x^2+20 x +5k \geq 0$, we see $k=13$
(this already shows the gonality is $13$).
Furthermore $16 \geq \phi(A,A)$ forces
 $(b,c)$ to be either $(0,0)$ or $(0,1)$ or $(1,0)$.
 Now $b=c=0$ gives
$-4=x=5a$, impossible, while $b=0$, $c=1$ gives $-4=x=5a+4$, also
impossible. Thus we must have $b=1$ and $c=0$. Then
$-4=5a+1$, hence $a=-1$. Therefore $A=L-H$, and then one easily checks that $\mathcal{Z}$ is precisely
the pencil cut out by planes through $L$ -see
Corollary~\ref{ms}. Thus the $g^1_{13}$ is unique.

\section{Notation and terminology} \label{notation}
A linear system of degree $k$ and projective dimension $r$ on $C$ is traditionally denoted
with the symbol $g_k^r$, and a $g_k^1$ is called a pencil. The {\em gonality} of $C$, written $\gon (C)$, is the least positive integer
$k$ such that there exists a $g^1_k$ on $C$. Since a pencil of least degree is automatically
base point free, the gonality of $C$ is the least degree of a surjective morphism $C \ra \PP^1$.
One can further notice that a $g^1_k$ with $k= \gon (C)$ is complete: this means that, if
$Z$ is a divisor in the given pencil,  $h^0 (C,\coo_C(Z))=2$, so that the pencil is the complete linear
series $|Z|$ of effective divisors linearly equivalent to $Z$.

\bd \label{dms}
Assume $C \subset \PP^3$ is a non planar curve.
Given a line $L$, let $\pi_L:C \ra \PP^1$ be obtained projecting $C$ from $L$,
and let
$\mathcal{Z}(L)$
denote the $g^1_k$ corresponding to $\pi_L$. Note that $\mathcal{Z}(L)$ is obtained
from the pencil cut out on $C$ by planes through $L$ removing its base locus, which coincides
with the scheme theoretic intersection $C \cap L$. In particular,
$$
\deg (\pi_L)= \deg \mathcal{Z}(L) = \deg (C) - \deg (C.L)
$$
and $\mathcal{Z}(L)$ is complete if $\deg (C.L) \geq 2$.
We say that a $g^1_k$ on $C$ arises from a multisecant if it is of the form $\mathcal{Z}(L)$ for some line $L$.
We say the gonality of $C$ can be {\it computed
by multisecants}
if there exists a line $L$ such that $\mathcal{Z}(L)$ has degree $\gon(C)$.
\ed

\section{Existence of $4$-secant lines} \label{Cayley}
The following statement is probably classical and well known, but it seems hard to find a reference.
\bp \label{Cay}
Let $C$ be a smooth irreducible curve of degree $d \geq 10$ in $\PP^3$. Then $C$ has an $l$-secant line $L$ with $l \geq 4$.
In particular, the gonality of $C$ is at most $d-4$.
\ep
\begin{proof}
We have to prove there is a line $L$ meeting $C$ in a scheme of length at least $4$.
This is clear if $C$ is a plane
curve of degree $\geq 4$, or if $\deg (C) \geq 7$ and $C$ is contained in a quadric surface. Thus we may assume
$C$ is not contained in a quadric surface.
In this case, we will show the Cayley number of $4$-secants
$$
\mathcal{C}(d,g)= \frac{(d - 2)(d - 3)^2(d - 4)}{12}-   \frac{g(d^2  - 7d + 13 - g)}{2}
$$
is positive. The existence of $L$ then follows from intersection theory as explained in Le Barz \cite{LB} or in \cite{acgh}.
For fixed $d \geq 7$, the number $\mathcal{C}(d,g)$ is a decreasing function of $g$,
because the partial derivative with respect to $g$ is
$$g-\frac{d^2  - 7d + 13 }{2}$$
which is negative because
$$g \leq  d^2/4-d+1$$ when $C$ is not contained in a plane.

But $C$ is not even contained in a quadric surface, thus its genus is bounded above by $\ds \frac{1}{6} d(d-3) + 1 $,
and
$$
\mathcal{C}(d,g) \geq  \mathcal{C}(  d,   \frac{1}{6} d(d-3) + 1   ) = \frac{d(d-3)(d-6)(d-9)}{72}
$$
which is positive for $d \geq10$.

\end{proof}
\br
The result is sharp, because a smooth complete intersection of two cubic surfaces has degree $9$ and no $4$-secant line.
\er

\section{Non existence of $5$-secant lines} \label{no5section}
\bt \label{no5secants}
Let $C \subset \PP^3$ be a curve
contained in an irreducible surface $X$ of degree $s$.
Suppose $C$ is a Cartier divisor on $X$, and
$$
\begin{cases}
H^0 (\PP^3, \ideal{C} (s\!-\!2)\,)= 0  \\
H^1 (\PP^3, \ideal{C} (m)\,)= 0 \mbox{\ \ \ \ for $m=s\!-\!2, s\!-\!3, s\!-\!4$}
\end{cases}
$$
If $C$ is general in its linear system on $X$,
then $\deg (C.L) \leq 4$  for every line $L$ not contained in $X$, and $C$ has only finitely many $4$-secant
lines not contained in $X$.

In particular, if $X$ does not contain a line, then $C$ does not have an $l$-secant line for any $l \geq 5$.
\et
\begin{proof}
The statement is obvious if $s \leq 3$, thus assume $s \geq 4$.

The hypotheses imply $h^1 \coo(D)=0$ for
$D=C$, $C-H$, $C-2H$ because by Serre duality
$$
h^1 (\PP^3, \ideal{C} (m))=
h^1(X, \coo_X (mH-C))= h^1 (X, \coo_X (C+(s-4-m)H)).
$$

Similarly, $H^2 (\coo_X(C-nH))$ is dual
to $$H^0 (\coo_X((s-4+n)H-C)) =H^0 (X,\ideal{C,X} (s-4+n))$$ which by assumption is zero for $n \leq 2$.
Thus we see that $h^0 \coo_X(D) = \chi \coo_X(D)$ for  $D=C$, $C-H$, $C-2H$.

Let $L$ be a line not contained in $X$, and let $V$ be the scheme theoretic intersection of $X$ and $L$.
Then $V$ has degree $s$, and there is an exact sequence
$$
\exact{\coo_X(-2H)}{\coo_X(-H)^{\oplus 2}}{\ideal{V,X}}.
$$
Twisting by $\coo_X(C)$ and taking cohomology we see
$$
h^0 (\ideal{V}(C)) = 2 h^0 (\coo_X(C-H))-h^0( \coo_X(C-2H))
$$

Therefore
\begin{eqnarray*}
&&h^0 (\coo_X(C)) -h^0 (\ideal{V}(C))=
h^0 (\coo_X(C))-2 h^0 (\coo_X(C\!-\!H))+h^0( \coo_X(C\!-\!2H)) =  \\
&&=\chi (\coo_X(C))-2 \chi (\coo_X(C-H))+\chi( \coo_X(C-2H)) =s
\end{eqnarray*}

This shows points of $V$ impose independent conditions on the linear system
$|C|$. It follows that the family of curves in $|C|$ meeting $L$ in a scheme
of length $l \leq s$ has codimension $l$ in $|C|$. This implies the statement because
$L$ varies in a $4$-dimensional family.
\end{proof}

\bc \label{no5s}
Let $C \subset \PP^3$ be an ACM curve. Suppose that $C$ is
contained in a smooth surface $X \subset \PP^3$ of degree $s=s_C$,
and that $C$ is general in its linear system on $X$.
Then $\deg (C.L) \leq 4$  for any line $L$ not contained in $X$.

In particular, if $X$ does not contain a line, then $C$ does not have an $l$-secant line for any $l \geq 5$.
\ec
\begin{proof}
The statement follows from Theorem \ref{no5secants} because $C$ is ACM precisely when $H^1 (\PP^3, \ideal{C} (m))=0$ for every $m$.
\end{proof}


\section{Gonality of curves on a smooth surface: Lazarsfeld's method } \label{laz}
In this section we explain a construction due to Lazarsfeld \cite{LK3,L} that will be crucial in proving that every pencil
of minimal degree on a general ACM curve arises from a multisecant.

When a curve $C$ is contained in a smooth surface $X$,  we associate a rank two vector bundle on $X$ to a base point free $g^1_k$ on $C$ as follows.
The base point free $g^1_k$ is determined by a degree $k$ line bundle $\coo_C(Z)$ on $C$, and a surjective map of $\coo_C$-modules
$$
\beta: \coo_C^{\, \oplus 2} \ra \coo_C (Z)
$$
(note that, since $k \geq 1$, the map $H^0 (\beta): H^0 (\coo_C^{\oplus 2}) \ra H^0 (\coo_C (Z))$ is
injective).

\bd \label{assoc}
Suppose $C$ is an integral curve on the smooth projective surface $X$, and $\mathcal{Z}$ is a base point free pencil
on $C$ defined by $\beta: \coo_C^{\, \oplus 2} \ra \coo_C (Z).$
Let $\alpha: \coo_X^{\, \oplus 2} \ra \coo_C (Z)$ denote the map obtained composing $\beta$ with the natural surjection
$\coo_X^{\, \oplus 2} \ra \coo_C^{\, \oplus 2}$.
Then the kernel $\mathcal{E}$ of $\alpha$ is called the {\em bundle associated} to the pencil $\mathcal{Z}$.
\ed

\bp \label{lazass}
Let $\mathcal{E}$ be the bundle associated to a  pencil of degree $k$ on $C$ as in the previous definition. Then
\begin{enumerate}
    \item[a)] $\mathcal{E}$ is a rank two vector bundle on $X$.
    \item[b)] $H^0 (\cae)=0$.
    \item[c)] $c_1 (\mathcal{E})=\coo_X (-C) \quad \mbox{and} \quad c_2 (\mathcal{E})= \deg (Z)$, so that
    $$
\Delta (\mathcal{E})\stackrel{def}{=} c_1^2(\mathcal{E})-4 c_2 (\mathcal{E})=C^2-4k
    $$
    \noindent(here we consider the first Chern
    class as an element of $A^1(X) \cong Pic(X)$, while we view the $c_1^2$ and $c_2$ as integers,
    via the degree map for zero cycles).
%
\end{enumerate}
\ep
\begin{proof}
By definition of $\mathcal{E}$ there is an exact sequence:
$$
\exact{\mathcal{E}}{\coo_X^{\, \oplus 2} }{\coo_C (Z)}
$$
Since $\coo_C$ has rank zero and projective dimension $1$ as an $\coo_X$-module, $\mathcal{E}$ is a rank two vector
bundle on $X$, whose Chern classes can be computed from the above sequence.
If $H^0 (\mathcal{E})$ were not zero, then $H^0 (\alpha): H^0 (\coo_C^{\, \oplus 2}) \ra H^0 (\coo_C (Z))$ would not be injective,
so $\alpha$ would induce a surjective map $\coo_C \ra \coo_C(Z)$, contradicting $\deg Z=k \geq 1$.
\end{proof}

We recall the definition of Bogomolov instability for rank two vector bundles on a surface, and Bogomolov's Theorem
which gives a numerical condition for instability.
\bd \label{unstable}
Let $\mathcal{E}$ be a rank two vector bundle on $X$.
One says that $\mathcal{E}$ is {\em Bogomolov unstable}
if there exist a finite subscheme $W \subset X$ (possibly empty) and divisors $A$ and $B$ on $X$
sitting in an exact sequence
\begin{equation}\label{bog}
\exact{\coo_X(A)}{\mathcal{E}}{\ideal{W}\otimes \coo_X(B) }.
\end{equation}
where $(A-B)^2 >0$ and $(A-B).H >0$ for some (hence every) ample divisor $H$. We say $A$ is a {\em destabilizing
divisor} of $\mathcal{E}$. It is unique up to linear equivalence.
\ed

\bt[Bogomolov \cite{Bog}, cf. {\cite[7.3.3]{HL}} and {\cite[4.2]{L}}] \label{bogomolov}
Suppose the ground field $\mathbb{K}$ has characteristic zero. Let $\mathcal{E}$ be a rank two vector bundle on
the smooth projective surface $X$, and let $\Delta (\cae)= c_1(\cae)^2-4 c_2 (\cae)$.

If $\Delta (\cae) >0$, then $\mathcal{E}$ is Bogomolov unstable.
\et

Following Lazarsfeld's approach, we will show in Section \ref{liaison} that the bundle associated to a pencil computing
the gonality of a smooth ACM curve satisfies $\Delta (\cae) >0$, hence it is Bogomolov unstable,
and there is a destabilizing divisor $A$. To work effectively we will need the following
technical result that will be useful in two ways. First it immediately
implies that, when $-A=H$ (plane section) or $-A=H-L$ (plane section minus a line),
the given pencil arises from a multisecant; later on the inequalities $A^2 \geq 0$ and $A.H<0$ will
be used to exclude all other possibilities for $A$.

\bp \label{a2}
Suppose $X$ is a smooth projective surface, $C$ is an integral curve on $X$, and
$|Z|$ is a complete base point free pencil on $C$. Let $\mathcal{E}$ be the rank $2$ bundle on $X$ associated to $|Z|$.
Suppose there is an exact sequence
\begin{equation}\label{dbog}
0 \rrr \coo_X(A) \stackrel{h}{\rrr} \mathcal{E}\rrr \ideal{W}\otimes \coo_X(B) \rrr 0
\end{equation}
with $W$ zero dimensional and $B$ not effective. Then
the linear system $|-A|$ on $X$ contains two effective curves $D_1$ and $D_2$ with the following properties:
\begin{itemize}
  \item[a)] $D_1$ and $D_2$ meet properly in a $0$-dimensional scheme $V$ containing $W$.
  \item[b)] $D_1$ and $D_2$ meet $C$ properly, and,
  if $R$ is the base locus of the pencil cut out on $C$ by $C.D_1$ and $C.D_2$, then
  $$\coo_C (Z) \cong \coo_X (-A) \otimes \coo_C(-R)$$
  i.e. the pencil $|Z|$ is obtained by first restricting $D_1$ and $D_2$ to $C$ and then removing the base locus $R$.
  \item[c)] $R$ is the residual scheme to $W$ in $V$, that is, there is an exact sequence
  $$
  \exact{\coo_W}{\coo_V}{\coo_R}
  $$
\end{itemize}
In particular $h^0 \ideal{W} (-A) \geq 2$,
$A.H <0$ for every ample divisor $H$, and $A^2 \geq 0$.
\ep
\br \label{a2r}
The proposition applies if $\mathcal{E}$ is Bogomolov unstable  with destabilizing sequence
(\ref{dbog}). Indeed in this case, if $H$ is an ample divisor on $X$, then $(A-B).H >0$.
Since $c_1 (\mathcal{E}) =A+B= -C$ in $\mbox{Pic} (X)$, we compute
$$-2B.H=(A-B).H+C.H >0.$$
Therefore $B$ is not effective.
\er

\begin{proof}[Proof of \ref{a2}]
Dualizing
$\exact{\mathcal{E}}{\coo_X^{\, \oplus 2}}{\coo_C (Z)}$
we obtain an exact sequence
$$
\exact{\coo_X^{\, \oplus 2}}{\mathcal{E}(C)}{\coo_C (C-Z)}.
$$
We now look at the composite map $g:\coo_X^{\, \oplus 2} \ra \mathcal{E}(C) \ra \ideal{W} (-A)$.

This map is nonzero, otherwise $\coo_X^{\, \oplus 2}$ would map injectively into
the kernel of $ \mathcal{E}(C) \ra \ideal{W} (-A)$, which is $\coo_X (C+A)$, absurd. Hence
the image of $g$ has rank one, and has the form $\ideal{Y} (-A)$ for some proper
subscheme $Y \subset X$ containing $W$. Then $\ideal{Y}=\ideal{V}(-D)$ where $D$ is the
divisorial part of $Y$, and $V$ is zero dimensional. We obtain an exact sequence
$$
\exact{\mathrm{Ker} (g)}{\coo_X^{\, \oplus 2}}{\ideal{V}(-A-D)}.
$$
It follows  $\mathrm{Ker} (g)= \coo_X(A+D)$ and $-A-D$ is effective. A diagram chase shows there is an exact sequence
$$
0 \ra \coo_X(A+D) \ra \coo_X (C+A) \ra \coo_C(C-Z)
$$
from which we see there is an effective curve $C_0$ linearly equivalent to $C-D$ contained in $C$.
Since $C$ is irreducible, this implies either $D=C$ or
$D=0$.

Now $-A-D$ is effective, so, if we had $D=C$, then $B=-A-C$ would be effective,
contradicting the hypotheses.
Hence the only possibility is $D=0$.

Putting everything together we obtain a commutative diagram with exact rows:
$$\begin{CD}
0 @>>> \coo_X(A) @>g^{\vee}>> \coo_X^{\, \oplus 2} @>(s_1,s_2)>> \ideal{V}(-A) @>>> 0\\
@. @VVhV @| @VVV @.\\
0 @>>> \mathcal{E} @>>> \coo_X^{\, \oplus 2} @>>> \coo_C(Z) @>>> 0
\end{CD}$$
Now let $D_1$ and $D_2$ the divisors defined by the sections $s_1$ and $s_2$ of $\coo_X(-A)$.
The first row of the diagram shows $D_1$ and $D_2$  meet properly in the zero dimensional scheme
$V$, which contains $W$ by construction. The two sections remain independent in $H^0 (\coo_C (Z))$ because $H^0 (\mathcal{E})=0$.
Hence $D_1$ and $D_2$ meet $C$ properly, and $D_1.C$ and $D_2.C$ span a pencil on $C$.

By the snake lemma, the kernel of the vertical map $\ideal{V}(-A) \ra \coo_C (Z)$ is $\ideal{W}(B)=\ideal{W}(-A-C)$,
hence a diagram chase produces an exact sequence
$$
\exact{\coo_C(Z)}{\coo_X(-A) \otimes \coo_C}{\coo_V / \coo_W}
$$
which proves the rest of the statement.
\end{proof}

\bc \label{ms}
Assume $X \subset \PP^3$ is a smooth surface with plane section $H$,
containing a smooth irreducible curve $C$.
Suppose $C$ is not contained in a plane. Let $|Z|$ be a complete base point free pencil on $C$,
and let $\mathcal{E}$ be the bundle on $X$ associated to $|Z|$.

\begin{itemize}
\item[a)]
If there is an exact sequence
\begin{equation*}\label{A-H}
\exact{\coo_X(A)}{\mathcal{E}}{\ideal{W}(B) }
\end{equation*}
with $W$ zero dimensional and $A+H$ effective, then
there is a line $L$ such that $|Z|=\mathcal{Z}(L)$ is the pencil cut out on $C$ by planes
through $L$. Furthermore, if $X$ does not contain $L$,
then $A=-H$ and $W$ is the residual scheme to $C \cap L$  in $X \cap L$, while, if $X$ contains $L$, then $A=L-H$ and $W$ is empty.
\item[b)] Assume $C$ is linearly normal and $|Z|$ is the pencil cut out on $C$ by planes
through a line $L$ meeting $C$ in a scheme of length at least $2$. Then there exists an exact sequence as above
with $A=-H$ if $X$ does not contain $L$ and $A=L-H$ if $X$ contains $L$.

\end{itemize}
\ec
\begin{proof} We first prove $a)$.
The divisor $B$ is not effective, otherwise
$$
B+(A+H)= (-A-C)+(A+H)=H-C
$$
would be effective, which contradicts the assumption that $C$ is not contained in plane.

Thus we may apply Proposition \ref{a2} to the given exact sequence
to conclude the linear system $|-A|$ contains a pencil. By assumption $P=A+H$ is effective,
and therefore in order that $|-A|=|H-P|$ may contain a pencil it is necessary
that $P$ be empty or a line.

If $P$ is empty, by \ref{a2} the are two plane sections $D_1=H_1 \cap X$ and $D_2=H_2 \cap X$  of $X$ meeting
in a zero dimensional scheme $V$, hence the line $L=H_1 \cap H_2$ is not contained in $X$.
Proposition \ref{a2}b shows $|Z|$ is obtained removing from the pencil spanned by $C \cap H_1$ and $C \cap H_2$ its base locus
$C \cap L$, that is, $|Z|=\mathcal{Z}(L)$, and Proposition \ref{a2}c shows $W$ is the residual scheme to $C \cap L$  in $X \cap L$.

Finally, if $P$ is a line, then $D_1$ and $D_2$ belong to $|H-P|$, hence their intersection $V=D_1 \cap D_2$ is empty.
It follows from Proposition \ref{a2} that $|Z|=\mathcal{Z}(P)$ and that
and  $W$ is empty. This proves the first statement.

We now prove $b)$. By definition of $\mathcal{E}$ there is an exact sequence:
$$
\exact{\mathcal{E}}{\coo_X^{\, \oplus 2} }{\coo_C (Z)}
$$
Comparing this sequence with
$$\exact{\coo_C}{\coo_C(Z)}{\coo_{Z}}$$
we obtain
$$
\exact{\coo_X (-C)}{\mathcal{E}}{\ideal{Z,X}}.
$$
Now twist by $H$ and take cohomology to get a long exact sequence:
$$
0 \ra H^0 (\coo_X (H-C)) \ra H^0 (\mathcal{E} (H)) \ra H^0 (\ideal{Z,X}(H))
\ra H^1 (\coo_X (H-C)).
$$
Since $Z$ is contained in a plane,  $h^0 (\ideal{Z,X}(H))>0$, while $H^1 (\coo_X (H-C))= H^1 (\ideal{C}(H))=0$ because
$C$ is linearly normal. Hence  $\mathcal{E} (H)$ has a section, and after removing torsion in the cokernel if necessary
we find an exact sequence:
$$
\exact{\coo_X (P-H)}{\mathcal{E}}{\ideal{W}(H-P-C)} \
$$
with $W$ zero dimensional and $P$ effective. Now $b)$ follows from $a)$.

\end{proof}

\section{ACM curves} \label{liaison}
In this section
we show that, if $C$ is an ACM curve of degree $d$ having a pencil of minimal degree
$k \leq d-4$ on a smooth surface of degree $s=s_C$, then the bundle $\mathcal{E}$ associated to the given pencil
satisfies $\Delta ({\mathcal{E}}) >0$ (except for a small list of cases \ref{delta-pos}), hence, if the ground field has characteristic zero,
it is Bogomolov unstable. The proof is based on
the structure of the biliaison class of ACM curves which we now briefly recall. We also include some information about the
minimal link $\Gamma$ of a curve $C$, which we will need later.

Given a curve $C$ in $\pso$ its fundamental numerical invariants are, besides its degree $d_C$ and
its arithmetic genus $g(C)=1 -\chi(\coo_C)$:
\begin{itemize}
    \item its index of speciality
    $$
    e(C)=\mbox{Max} \{n: h^1 \coo_C (n) >0\};
    $$
    \item the minimal degree $s_C$ of a surface containing $C$;
    \item the integer
    $$t_C= \mbox{Min} \{n: h^0 (\ideal{C}(n))-h^0 (\coo_{\pso}(n-s_C))>0)\}$$
    If $C$ is integral or more generally if $C$ lies on an integral surface of degree $s_C$, the integer $t_C$ is
the smallest $n$ such that $C$ is contained in a complete intersection of two surfaces of degree $s_C$
and $n$.
\end{itemize}

When $C$ is ACM, all its basic numerical invariants can be computed from the Hilbert function. It is convenient to express the Hilbert function
through its second difference function, the so called {\em $h$-vector} $h_C$ of $C$ - cf. \cite[\S 1.4]{M} -
because $h_C$ is a finitely supported function.
Thus one defines
$$
h_C (n)=
h^0(\coo_C (n))-2 h^0(\coo_C(n-1))+h^0(\coo_C(n-2)).
$$
If $s=s(C)$ and $e=e(C)$, the function $h_C$ satisfies
\begin{equation} \label{admiss}
\begin{cases} h(n) = n+1 \;\;\;\mbox{if \ } 0 \leq n \leq s\!-\!1, \\
h(n) \geq h(n+1) \;\;\;\mbox{if \ } n \geq s\!-\!1, \\
h(e+2)>0 \quad \mbox{and} \quad h(n)=0 \quad \mbox{for $n \geq e+3$}.
\end{cases}
\end{equation}
Thus we may write $h$ as
 \begin{equation*}\label{hvettore}
    h_C = \{ 1,2, \ldots, s, h_C (s),  \ldots, h_C (e+2) \}.
\end{equation*}
with $s=h_C(s\!-\!1) \geq h_C(s) \geq h_C(s\!+\!1) \geq \ldots \geq h_C(e+2)$.

We say that a finitely supported function $h: \N \ra \N$ is an  {\em $h$-vector}
if it satisfies (\ref{admiss}) for some $s\geq 1$. Every $h$-vector arises
as the $h$-vector of an ACM curve in $\PP^3$ (cf. \cite[Theorem V.1.3 p. 111]{MDP} and Remark \ref{exist} below). It will be convenient
to allow the identically zero function among $h$-vectors, and think of it as the $h$-vector of the empty curve.

In terms of the $h$-vector, the fundamental invariants of $C$ are computed as follows:
\bp
Let $C$ be an ACM curve in $\pso$ with $h$-vector $h_C$. Then
\begin{enumerate}
    \item $\ds d_C=  \sum h_C(n)$;
    \item $\ds g(C)= 1+ \sum (n-1) h_C(n)$;
    \item $\ds e(C)+2= \mbox{\em Max} \{n: h_C (n) >0 \}$;
    \item $\ds s_C= \mbox{\em Min} \{n \geq 0: h_C(n) < n+1 \}$;
    \item
    $ \ds t_C= \mbox{\em Min} \{n \geq 0: h_C(n-1) > h_C(n) \}$.
\end{enumerate}
\ep
Coherently with the above formulas, for the empty curve we define $s=0$, $d=0$, $g=1$, $e= - \infty$.

\br \label{degreebound}
If $C$ is an ACM curve with $s_C=s$, then
$$
d_C = \sum h_C(n) \geq \sum_{n=0}^{s\!-\!1} (n+1)= \frac{1}{2} s (s\!+\!1).
$$
\er

The $h$-vectors of integral curves have a special form:
\bd[see\cite{MR}]
An $h$-vector is of {\em decreasing type}  if
$h(a)>h(a+1)$ implies that for each $n \geq a$ either $h(n)>h(n+1)$ or $h(n)=0$.
\ed

\br \label{ell}
We also recall that, by a Theorem of Ellingsrud's
the Hilbert scheme $A(h)$ of ACM curves in $\pso$ with a given $h$-vector is smooth and irreducible,
even when $h$ is not of decreasing type - see \cite{E} and \cite[p. 5, Corollaire 1.2 p. 134 and 1.7 p. 139]{MDP}.

Gruson and Peskine \cite{GP} (see also \cite{MR} and \cite{nollet})
showed that, if $C$ is an integral ACM curve, then $h_C$ is of decreasing type, and
conversely, if $h$ is an $h$-vector of decreasing type, then there exists a smooth irreducible ACM curve
$C$ with $h_C=h$.
Thus an $h$-vector $h$ is of decreasing type if and only if the general curve $C$ in
$A(h)$ is smooth and irreducible.
\er
If $C$ is not irreducible, it may happen that every pair of surfaces $X_1$ and $X_2$ containing $C$ of minimal degrees
$s_C$ and $t_C$ have a common component. Nollet \cite[Proposition 1.5]{nollet}
generalizes the result of Gruson and Peskine by showing that if $C$ is contained in a complete intersection
of type $(s_C,t_C)$, then $h_C$ is of decreasing type. We partially reproduce his argument in the following lemma:
\bl \label{decreasing} \

\begin{itemize}
    \item[i)]
    Suppose an ACM curve $D$ is contained in a complete intersection $Y$ of type $(s_D,t_D)$, and let $\Gamma$ be the
curve and linked to $D$ by $Y$. Then $$e(\Gamma)+3 <s_D.$$
    \item[ii)]
    Let $\Gamma$ be an ACM curve, and suppose $a \leq b$ are integers such that $a \geq e(\Gamma)+3$ and $b \geq e(\Gamma)+4$.
Then the $h$-vector of a curve $D$ linked to $\Gamma$ by a complete intersection of type $(a,b)$ is of decreasing type.
If $a \geq e(\Gamma)+4$, then $s_D=a$ and $t_D=b$. If $a=e(\Gamma)+3$, then $s_D=a$ and $t_D=b-1$.
\end{itemize}
\el
\begin{proof}
If $\Gamma$ and $D$ are linked by a complete intersection $Y$ of type $(a,b)$, then by \cite[5.2.19]{M}
$$
h_{\Gamma} (n)= h_Y (n)-h_D (a+b-2-n) = h_Y (a+b-2-n)-h_D(a+b-2-n).
$$

Suppose first $a=s_D$ and $b=t_D$. Then
$$
h_{\Gamma} (s_D-1)= h_Y (t_Y-1)-h_D(t_D-1)=s_Y-s_D=0.
$$
Therefore $e(\Gamma)+3 \leq s_D-1$.

Next suppose $b \geq a \geq e(\Gamma)+4$.
Then $s_D \leq a$ because $D \subseteq Y$, and
$$
h_D(b-1)=h_Y(a-1)-h_\Gamma (a-1) = h_Y (a-1)=a
$$
while
$$h_D (b)= h_Y(a-2)-h_\Gamma (a-2) \leq h_Y (a-2)=a-1$$
hence $s_D=a$ and $t_D=b$.

If $a=e(\Gamma)+3$ and $b \geq e(\Gamma)+4$, then
a similar calculation shows $h_D (b-2)=a$ and $h_D (b-1) <a$, so that
$s_D=a$ and $t_D=b-1$.

It remains to show $h_D$ is of decreasing type. Let $u=s(\Gamma)$.  Then
$ u \leq e(\Gamma)+3 \leq a$ and
$$h_\Gamma(n)=h_Y (n)=n+1 \quad \mbox{for $n \leq u-1$,}$$
hence $h_D (n)=0$ for $n \geq a+b-1-u$.

Since $h_\Gamma (n) \geq h_\Gamma(n+1)$ for $n \geq u-1$, we see that for $b-1 \leq m \leq a+b-2-u$
\begin{eqnarray*}
h_D(m)-h_D(m+1) &=& h_Y(m)-h_Y(m+1)-h_\Gamma (a+b-2-m) + \\+ \, h_\Gamma(a+b-1-m)
&=&1- \partial h_\Gamma (a+b-1-m) \geq 1
\end{eqnarray*}
which shows that $h_D$ is of decreasing type.
\end{proof}

Fix a smooth surface $X \subset \PP^3$ of degree $s$.
Two curves $C$ and $D$ on $X$ are said to be {\em biliaison equivalent}
if $C$ is linearly equivalent to $D+nH$ for some integer $n$.
\bd
A curve $C$ on a surface $X$ is {\em minimal on $X$} if $C-H$ is not effective.
\ed

\bp \label{minimal}
A curve $C$ is minimal on a smooth surface $X$ if and only if $e(C)+3< \deg (X)$
\ep
\begin{proof}
To say $C$ is minimal is equivalent to saying $h^0 (\coo_X (C-H))=0$. By duality on $X$ this is the same
as $h^2 (\ideal{C} (s\!-\!3))=0$, where $s=\deg (X)$. On the other hand, $h^2 (\ideal{C} (s\!-\!3))=h^1 (\coo_{C} (s\!-\!3))$,
so the condition says $s-3>e(C)$,i.e. $e(C)+3 <s$.
\end{proof}

\bd
We say that an $h$-vector is {\em $s$-minimal} if the corresponding curve
satisfies $e +3<s$.
We say that an $h$-vector is {\em $s$-basic} if it is the $h$-vector of an integral curve
$C$ satisfying $s_C=t_C=s$. Thus the $s$-basic $h$-vectors are those $h$-vectors of decreasing type that begin with a string
$$
\{ 1,2, \ldots , s\!-\!1,s,m \}
$$
with $m =h(s) \leq s\!-\!1$.
\ed

Figure \ref{table1} at the end of the paper lists $s$-basic $h$-vectors for $s=4$ and $s=5$.
\bp \label{acm}
Suppose $C$ is an ACM curve contained in a smooth surface $X$ of degree $s_C$.
Let $s=s_C$, $t=t_C$ and $e=e(C)$.
Then $ e+3 \geq t \geq s$ and
\begin{itemize}
    \item[a)] $h_C$ is of decreasing type;
    \item[b)] if $\Gamma \in |tH-C|$, then $e(\Gamma)+3 <s$ and $\Gamma$ is minimal on $X$;
    \item[c)] $C-mH$ is effective if and only if $m \leq e+4-s$;
    \item[d)] if $C_1 \in |C-(t\!-\!s)H|$, $h_{C_1}$ is $s$-basic;
    \item[e)] if $C_2 \in |C-(t\!-\!s\!+\!1)H|$, $h_{C_2}$ is of decreasing type.
\end{itemize}

There is a one to one correspondence $h_\Gamma \mapsto h_{C_1}$ mapping
$s$-minimal $h$-vectors to $s$-basic $h$-vectors.
\ep
\begin{proof}
Since $C$ is ACM, the ideal sheaf $\ideal{C,\PP^3}$ is $(e+3)$-regular, hence
$e+3 \geq t$. By definition of $t$, we have $t \geq s$, and $C$ is contained
in a surface $F$ of degree $t$ that does not contain $X$. Therefore $C$ is contained
in the complete intersection $X \cap F$ of type $(s,t)$. Let $\Gamma \in |tH-C|$ be the curve
linked to $C$ by $X \cap F$: then $e(\Gamma_0)+3 <s$ and $\Gamma$ is minimal (by either Lemma \ref{decreasing} or
by definition of $t$).

Note that $C$ (respectively $C_1$, resp.  $C_2$) is linked to a curve in the linear system $|\Gamma|$ by a complete intersection
of type $(s,t)$ (resp. $(s,s)$, resp.  $(s\!-\!1,s)$)
By Lemma \ref{decreasing} the $h$-vectors of $C$, $C_1$ and $C_2$ are of decreasing type, and $h_{C_1}$ is $s$-basic.

\end{proof}

There is a unique $1$-basic $h$-vector, namely $h_0=\{1\}$, the $h$-vector of a line.
Every $(s\!-\!1)$-basic $h$-vector gives rise to two $s$-basic $h$ vectors by performing
a type $A$ or type $B$ transformation, where
\begin{enumerate}
    \item A type $A=A_s$ transformation consists of inserting an $s$ to an $(s\!-\!1)$-basic h-vector
$h=\{ 1,2, \ldots , s\!-\!1,m, .... \}$ to transform it into the $s$-basic vector
$h'=\{ 1,2, \ldots , s\!-\!1,s,m .... \}$. Geometrically, if $h$ is the $h$-vector of a curve $C$ on a surface
$X$ of degree $s$, $h'$ is the $h$-vector of the effective divisor $C+H$ on $X$.
    \item  A type $B$ transformation  consists of inserting a string $s,s\!-\!1$ to an $(s\!-\!1)$-basic h-vector
$h=\{ 1,2, \ldots , s\!-\!1,m, .... \}$ to transform it into the $s$-basic vector
$h''=\{ 1,2, \ldots , s\!-\!1,s,s\!-\!1,m .... \}$. Geometrically, this operation breaks into two steps:
suppose $h$ is the $h$-vector of a curve $C$ on a surface
$X_1$ of degree $s\!-\!1$. Let $C_1=C+H$ be obtained by adding to $C$ a plane section of $X_1$, then pick a surface $X_2$ of degree $s$ containing $C_1$, and finally let $C_2=C_1+H$ be obtained by adding to $C_1$ a plane section of $X_2$. Then $h''$ is the $h$-vector of $C_2$.
\end{enumerate}
Conversely, any $s$-basic $h$-vector with $m=h(s) \leq s-2$ arises from a type $A$ transformation of
an $(s\!-\!1)$-basic $h$-vector, while any $s$-basic $h$-vector with $m=h(s)= s\!-\!1$ arises from a type $B$ transformation of
an $(s\!-\!1)$-basic $h$-vector.  In particular, the number of $s$-basic $h$-vectors is $2^{s\!-\!1}$ (cf. Figure \ref{table1}).

\bp \label{delta-pos}
Let $C$ be an integral ACM curve in $\PP^3$ with $s_C \geq 4$. Suppose $C$ is contained in a smooth
surface $X$ of degree $s=s(C)$. Suppose $C$ has a base point free pencil of degree $k$, and let $\mathcal{E}$ be the bundle on $X$ associated to such a pencil.

Then:
\begin{enumerate}
  \item[a)] if $k \leq d-5$, then $\Delta (\mathcal{E})>0$ unless
\begin{itemize}
    \item $s=4$ and $(d,g)=(10,11)$, or
    \item $s=5$ and $(d,g)=(15,26),(16,30)$, or
    \item $s=6$ and $(d,g)=(21,50),(22,55),(23,60)$, or
    \item $s=7$ and $(d,g)=(28,85),(29,91)$, or
    \item  $s=8$ and $(d,g)=(36,133)$.
\end{itemize}

  \item[b)] if $k=d-4$, then $\Delta (\mathcal{E})>0$ unless
  \begin{itemize}
    \item  $s=4$ and $(d,g)=(10,11), (11,14),(12,17)$, or
    \item  $s=5$ and $(d,g)=(15,26),(16,30),(17,34),(18,38)$, or
    \item  $s=6$ and $(d,g)=(21,50),(22,55),(23,60),(24,65)$, or
    \item  $s=7$ and $(d,g)=(28,85),(29,91), (30,97)$, or
    \item  $s=8$ and $(d,g)=(36,133),(37,140)$.
  \end{itemize}

\end{enumerate}

\ep
\begin{proof}
We can compute $\Delta (\mathcal{E})$ in terms of $d=d_C$ and $g=g(C)$:
$$
\Delta (\mathcal{E})= C^2-4k = 2g-2-(s-4)d-4k= \delta_s (d,g) +4 (d-k)
$$
where we have set
$$\delta_{s} (C)= \delta_s (d,g) =2g-2-ds.$$

One can easily verify that
\begin{enumerate}
    \item
Let $C \subseteq X_s$ be a curve on a surface $X$ of degree $s$ in $\PP^3$, and consider the divisor $C+H$ on $X_s$.
Then
$$
\delta_s (C+H)-\delta_s (C)= 2d-3s.
$$
In particular, if $d \geq \frac{1}{2}s(s\!+\!1)$ and $s \geq 3$, $\delta_s(C+H) > \delta_s (C)$.
   \item
Suppose $C \subseteq X_{s\!+\!1}$ is a curve on a surface $X$ of degree $s\!+\!1$ in $\PP^3$, and consider the divisor $C+H$ on $X_{s\!+\!1}$
Then
$$
\delta_{s\!+\!1}(C+H)-\delta_{s} (C)= d-3(s\!+\!1).
$$
In particular, if $d \geq \frac{1}{2}s(s\!+\!1)$ and $s \geq 6$, $\delta_{s\!+\!1}(C+H) \geq \delta_s (C)$, and
the inequality is strict unless $s=6$ and $d=21$.
\end{enumerate}

To prove the proposition, we have seen that $\Delta  (\mathcal{E})$ can be computed in terms of $d,g,s,k$, which depend
only on the $h$-vector and the choice of $s,k$. Therefore, using the two remarks $(1)$, $(2)$ just made and using biliaisons
on each surface to reduce to $s$-basic $h$-vectors, and using the transformations of type $A$ and $B$
mentioned before the statement, it would be sufficient to prove that $\Delta >0$ for all $s$-basic $h$-vectors with $s=4$.
Unfortunately this is not so, as $\Delta \leq 0$ for the first three $4$-basic $h$-vectors (see Figure \ref{table1}).
Still the two remarks show that $\Delta$ becomes positive using the transformations of type $A$ and $B$,
with the only exceptions listed in the statement. Figure \ref{table1} displays all $h$-vectors for which
$\Delta \leq 0$ for $k=d-4$ and $k=d-5$.
\end{proof}

\section{General ACM curves} \label{gacm}

In this section we give a description of a general ACM curve in the case when the $h$-vector is not of decreasing type,
thus generalizing the work of Gruson and Peskine \cite{GP}. We show (\ref{dthm}) that it is a union of smooth ACM subcurves that
are determined by the gaps in the associated biliaison type $\lambda$ (defined below). As a corollary we show the existence of
multisecant lines for ACM curves having certain particular behavior of the $h$-vector.

%

%
\bd\label{defbil}
Let $C_0$ and $C$ be two curves in $\PP^3$.
\begin{itemize}
    \item[a)] Following \cite{MDP} we say
that $C$ is obtained by an {\em elementary  biliaison} of height
$h$ from $C_0$ if there exists a surface $X$ in ${\mathbb
P}^3$ containing $C_0$ and $C$ so that
$\ideal{C,X} \cong \ideal{C_0,X}(-h)$. In the language of generalized divisors \cite{GD}
this means
$C$ is linearly equivalent to $C_0 + hH$ on $X$, where $H$
denotes the plane section.
    \item[b)] As a particular case, we say $C$ is obtained by a {\em trivial} biliaison
of height $h$ if $\ideal{C,X} = \ideal{C_0,X} \ideal{Y,X}$
where $Y$ is a complete intersection of $X$ and a surface of degree $h$. If $Y$ meets $C_0$ properly,
this means $C$ is the union of $C_0$ and $Y$.
    \item[c)] By a {\em special biliaison of degree $k$} we mean an elementary biliaison  of height one
     $C \sim C_0 + H$ on a surface of degree $ k \geq e(C_0)+4$. The condition $k \geq e(C_0)+4$ guarantees  $s_C=s_{C_0}+1$
and $k=e(C)+3$ by \cite[p. 68]{MDP}.
\end{itemize}
\ed

\bp[Lazarsfeld-Rao property] \label{LR}
Suppose $C$ is an ACM curve with index of speciality $e$.
Then $C$ can be obtained by a special biliaison
of degree $k=e+3$ from some ACM curve $C_0$ satisfying $s_{C_0}=s_C\!-\!1$.
\ep
\begin{proof}
One knows - see for example \cite{strano} - that an ACM curve $C$ with index of speciality $e$ can be obtained
by an elementary biliaison of height $1$ on a surface $X$ of degree $e+3$ from an ACM curve $C_0$
satisfying $$s_{C_0}=s_C\!-\!1 \quad \mbox{and} \quad
 e(C_0) < e(C).$$
Since $ \deg (X)=e+3 \geq e(C_0)+4$, this is a special biliaison.
\end{proof}

\br
When $s_C=1$, the curve $C_0$ above is the empty curve, which is therefore convenient to allow among ACM curves.
\er


\bc \label{btype}
Let $C$ be an ACM curve. Then there exist positive integers $k_1 <  k_2 <  \cdots <  k_u$
such that $C$ is obtained from the empty curve by a chain of $u$ special biliaisons of degrees $k_1, \ldots, k_u$.
The sequence $\lambda_C=(k_1, k_2, \ldots, k_u)$ is uniquely determined by $C$,
and we will call it the {\em biliaison type} of $C$. Furthermore:
\begin{enumerate}
    \item $\ds d_C=  \sum_{i=1}^{u} k_i$
    \item $\ds s_C= u$
    \item $\ds  g(C)= 1+\sum_{i=1}^{u} \frac{k_i( k_i-3)}{2} +
    \sum_{i=1}^{u} (s_C-i)k_i $
    \item $ \ds t_C-s_C+1= k_1$
    \item $\ds e(C)+3 = k_u$
\end{enumerate}
\ec

\bex \label{ds}
If $C \subset \PP^3$ is ACM , then $\ds d_C \geq \frac{1}{2}s_C(s_C+1)$,
with equality if and only if $\lambda_C=(1,  2,  3, \cdots , s_C-1,  s_C)$.
\eex

\br
The biliaison type $\lambda_C$ was introduced from a different point of view in
\cite{Gr}, and it essentially the same thing as the numerical character $\{n_j\}$ of Gruson and Peskine \cite{GP}: the precise
relationship, if $s=s_C$,  is
$$
n_j-j=k_{s-j} \;\;\;\;\;\; \mbox{for $j=0, \ldots, s\!-\!1$}.
$$

The biliaison type (hence the numerical character) is equivalent to the $h$-vector of $C$.
Indeed, $h_C$ can be recovered from
$\lambda_C$ because one knows how $h_C$ vector varies in an elementary biliaison, while $\lambda_C$ can be computed out of $h_C$ via the formula
$$
k_i = \#\{n: h_C (n) \geq s_C+1-i\}.
$$
One can visualize $h_C$ and $\lambda_C$ as follows. In the first quadrant of the $(x,y)$ plane,
draw a dot at $(n,p)$ if $n$ and $p$ are integers satisfying  $1 \leq p \leq h(n)$. Then $h(n)$ is the number
of dots on the vertical line $x=n$, while $k_i$ is the number of dots on the horizontal line
$y=s-i+1$. In particular,  $k_1=t_C-s_C+1 $ is the number of dots on the top horizontal line $y=s$,
and $k_s= e(C)+3$ is the number of dots on the bottom line $y=1$.
\er

\begin{figure}[htb]
\caption{Biliaison type and $h$-vector} \label{fig}
\begin{center}
\includegraphics[scale=1]{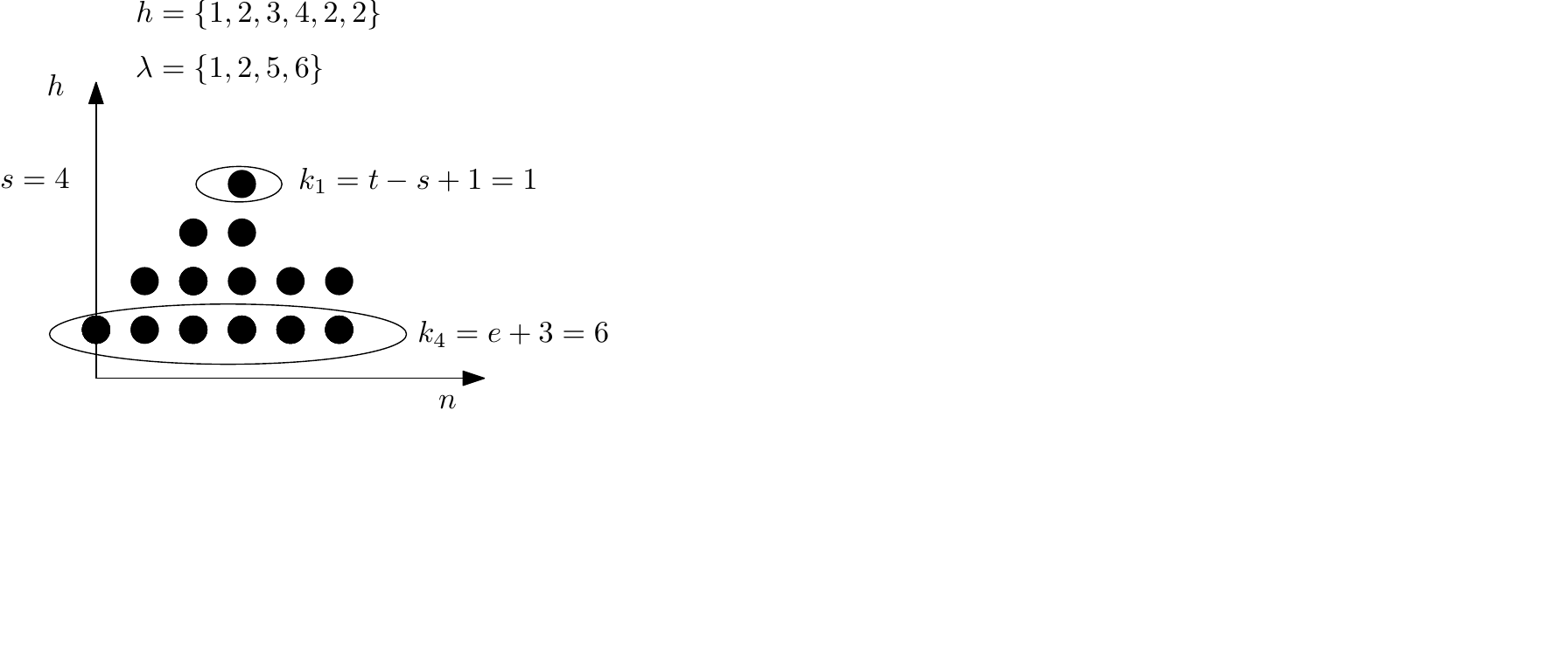}
\end{center}
\end{figure}

\br \label{exist}
The statement that every $h$-vector arises
as the $h$-vector of an ACM curve in $\PP^3$  is equivalent to the statement
 that every finite, strictly increasing sequence
of positive integers $\lambda=(k_1, \ldots, k_u)$  occurs as $\lambda_C$ for some ACM curve $C \subset \PP^3$. We can see
this by induction on $u$. When $u=1$, $\lambda=(k)$ is the biliaison type of a plane curve of degree $k$. If $u>1$,
by induction there is an ACM curve $C_0$ with $\lambda_{C_0}=(k_1, \ldots, k_{u-1})$. Now $s_{C_0}\leq e(C_0)+3=k_{u-1}<k_u$.
Therefore we can find a surface $X$ of degree $k_u$ containing $C_0$, and define $C$ to be a curve  obtained from $C_0$
by a biliaison of height one on $X$. Since $e(C_0)+3< k_u$, the biliaison is special, hence
$\lambda_C$ equals the given $\lambda$.
We give a refined version of this construction in Theorem \ref{dthm}.
   \er

\bd We say that a sequence $\lambda=(k_1,  k_2, ,  \ldots, k_u)$ {\em has a gap} at $i$
if $k_{i+1}-k_{i}\geq 3$.
\ed

For example, the sequence $\lambda_C$ of Figure \ref{fig} has a gap at $i=2$.
Note that
$\lambda$ has no gaps if and only if its corresponding $h$-vector is of decreasing type.
Gruson and Peskine \cite{GP} showed that, if $C$ is an integral ACM curve, then $\lambda_C$ has no gaps, and
conversely, if $\lambda$ has no gap, then there exists a smooth irreducible ACM curve $C$ with $\lambda_C=\lambda$.
Theorem \ref{dthm} below generalizes this statement to the case $\lambda$ has gaps; even Proposition \ref{gap}, which
says that $C$ is the union of two ACM subcurves whenever $\lambda_C$ has a gap, seems not to have been noticed before.
Here the new remark with respect to \cite{GP} is that the subcurves are themselves ACM. Note that if $B$ and $D$ are
two curves on a surface $X$ and $C$ is linearly equivalent to $B+D$, then it is possible that
$B$ and $C$ are ACM without $D$ being ACM (take for example $B$ to be a line, $D$ the disjoint
union of two lines each meeting $B$), or that $B$ and $D$ are ACM without $C$
being ACM (for example when $B$ and $D$ are two disjoint ACM curves).

So we need to impose some condition if we want the union of two ACM curves to be ACM. The condition that
suits our needs is that $I_D/I_C$ should be isomorphic to $R_B$ up to a twist.
This condition is satisfied when $C$ is obtained from $B$ by a trivial biliaison (see the examples below),
and also when $C$ is obtained from $B$ by a chain of elementary biliaison "that are trivial on $B$"
(Lemma \ref{commute} below).
Here are some preliminary examples.
\bex
If $C$ is obtained from a curve $B$ by a trivial biliaison of height $h$ on a surface $X$,
"adding" to $C$ the complete intersection $Y$  of $X$ with a surface of degree $h$, then
$$ I_Y/I_C \cong \frac{I_Y/I_X}{I_C/I_X} \cong
\frac{H^0_* (\ideal{Y,X}) }{H^0_* (\ideal{C,X})}
\cong
\frac{H^0_* (\coo_X(-h)) }{H^0_* (\ideal{B,X} (-h) )}
\cong R_B (-h)
$$
\eex
Here are more examples
\bex
Let $D \subset \pso$ be a curve, and $L$ a line not contained in $D$. Set $C=D \cup L$, and let $f$ be the
degree of the scheme theoretic intersection $D \cap L$. Then
$$\ideal{D,C} \cong \ideal{D \cap L,L} \cong \coo_L(-f).$$
Suppose now $D$ is ACM. Taking cohomology we obtain an sequence:
$$
0 \ra I_C \ra I_D \ra H^0_* (\ideal{D,C}) \cong R_L (-f) \ra H^1_* (\ideal{C}) \ra  H^1_* (\ideal{D})=0
$$
Hence $C=D \cup L$ is ACM if and only if $I_D$ maps onto $R_L (-f)$, giving an isomorphism
$I_D/I_C \cong R_L(-f)$.

By the same argument, if $B$ and $D$ are two ACM curves meeting properly and $\ideal{B \cap D, B} \cong \coo_{B} (-f)$,
then $C=B \cup D$ is ACM if and only if $I_D/I_C \cong R_B(-f)$.

From another point of view, suppose $B$ and $D$ are two ACM curves contained in a {\em smooth} surface $X$, and let $C=B+D$.
Then $$\coo_B (-D)=_{def}\coo_X(-D) \otimes \coo_B \cong \ideal{D,C}.$$
If  $\coo_B(-D) \cong \coo_B (-f)$, then $C$ is ACM if and only if $I_D/I_C \cong R_B(-f)$.
\eex

\bex
If $C=B \cup D$ and $B$, $C$ and $D$ are ACM, it may happen that $I_D/I_C$ is not isomorphic to
$R_B(-f)$ for any $f$. For example, let $B$ a twisted cubic curve, let $C$ be a complete intersection
containing $B$, and let $D$ the curve linked to $B$ by $C$. Then
$$\ideal{D,C} = \mathcal{H}om_{\coo_C} (\coo_B, \coo_C) \cong \omega_B (-e)$$
where $e=e(C)$. Hence $I_D/I_C$ is up to a twist $H^0_* (\omega_B)$, which requires two generators as an $R_B$
module.
\eex

\bex
If $B$ and $D$ have a common component, the conditions
that $C$ contains  $B$ and $D$, and
$I_D/I_C \cong R_B(-f)$, do not determine $C$.
For example, if $B=D=L$ is a line, $H$ is a plane containing $L$, and $C_H$ is the double structure on $L$
contained in $H$, then $\ideal{L,C_H} \cong \coo_L (-1)$ for every $H$.
\eex

The condition $I_D/I_C \cong R_B(-f)$ implies that $C$ is obtained by a "generalized liaison addition" of $B$ and $D$
in the sense of \cite{GM}. The following proposition is essentially a special case of Theorem 1.3 of \cite{GM}.
\bp \label{gen-add}
Suppose that $C$ contains two subcurves $B$ and $D$, and that for some integer $f$ there is an isomorphism of $R_C$-modules:
\begin{equation}\label{genadd}
I_D/I_C \cong R_B(-f).
\end{equation}
 Then
\begin{enumerate}
    \item[a)] There is a surface $S$ of degree $f$ containing $D$ but not $C$, and the curve $D$ is the scheme theoretic
    intersection of $C$ and $S$. In particular, $f \geq s_D$.
    \item[b)] The degrees and genera of $B$, $C$ and $D$ are related by the formulas
$$
    d_C=d_B+d_D  \;\;\;\;\;\;\;\; g(C)=g(B)+g(D)+f d_B -1
$$

     If $B$ and $D$ have no common component, then $C$ is the scheme theoretic union of $B$ and $D$,
    $\ideal{B \cap D, B} \cong \coo_{B} (-f)$, and $B.D=f d_B$.

    If $C$ is contained in a smooth surface $X$, then $C=B+D$ on $X$, and $\coo_X (D) \otimes \coo_B \cong \coo_B(f)$. In particular, $B.D=f d_B$ (here $B.D$ denotes the intersection product on $X$, and the formula holds
    even when $B$ and $D$ have common components).

    \item[c)] Suppose $D$ is ACM. Then $B$ is ACM if and only if $C$ is ACM, in which case
    $$
h_C (n)= h_B(n-f)+h_D(n)
    $$
    \item[d)] Suppose $B$, $C$ and $D$ are ACM and $f=s_D$.
    If $Max \{\lambda_B\} < Min \{\lambda_D\}$ then
    $$
    \lambda_C= \lambda_B \cup \lambda_D.
    $$
\end{enumerate}
\ep
\bex
Figure \ref{fig} shows the $h$-vector of a curve which is the union of a twisted cubic curve $B$
and a divisor of type $(6,5)$ on a smooth quadric surface as in Example \ref{example1}
In this case, $h_B =\{1,2\}$, $h_D=\{1,2,2,2,2\}$ and $f=s_D=2$. We see the $h$-vector
of $B$ sits on top of that of $D$. The biliaison types are
$\lambda_B= \{1,2 \}$ and $\lambda_D=\{ 5, 6\}$.
\eex
\begin{proof}[Proof of \ref{gen-add}]
The hypothesis $I_D/I_C \cong R_B(-f)$ is equivalent to the following; there is a  form
$F \in H^0 (\PP^3, \coo (f))$ such that the following sequence is exact:
$$
0 \ra I_B/I_C(-f) \ra R_C (-f) \stackrel{F}{\ra}  R_C \ra R_D \ra 0
$$
In particular, $I_D= I_C+I_S$ where $S$ is the surface of equation $F=0$,
hence $D$ is the scheme theoretic union of $C$ and $S$.
Sheafifying the exact sequence
$$
\exact{I_{B} (-f)}{I_{C}}{I_{D}/(F)}
$$
we obtain another exact sequence
$$
\lexact{H^1_* (\ideal{B}) (-f)}{H^1_* (\ideal{C})}{H^1_*(\ideal{D})}
$$
It follows that, if $D$ is ACM, then $H^1_* (\ideal{B}) (-f) \cong H^1_* (\ideal{C})$. In particular,
if $D$ is ACM, then $B$ is ACM if and only if $C$ is ACM  (generalizing the case of elementary biliaison).

If $B$ and $D$ are ACM, the relation between the $h$-vectors follows immediately from the exact sequence
$\exact{R_B(-f)}{R_C}{R_D}$.

The relation between the degrees and genera follows computing the Euler characteristics of the
two sides of  $\ideal{D,C} \cong \coo_B(-f)$.

Suppose $B$ and $D$ have no common components. The kernel of the natural surjective map
$$
\coo_B (-f) \cong \ideal{D,C}  \ra \ideal{B\cap D,B}
$$
is supported on $D$ and is a subsheaf of $\coo_B$. Since $B$ is locally Cohen-Macaulay and has no component
in common with $D$, the kernel is zero, hence $\coo_B (-f) \cong  \ideal{B\cap D,B}$.

Suppose $C$ is contained in a smooth surface $X$. Since $D \subseteq C$, there is an effective divisor $A$ on $X$
such that $C=A+D$. Then $$\coo_B (-f)\cong \ideal{D,C} \cong \coo_X (-D) \otimes \coo_A$$
from which we deduce $A=B$ and $\coo_B (f) \cong \coo_X (D) \otimes \coo_B $, hence $B.D=f d_B$.

We deduce  $d)$ from   $c)$. By assumption
$$e(B)+3=Max \{\lambda_B\} < Min \{\lambda_D\} = t_D-s_D+1.
$$
On the other hand, $h_D(n)= s_D$ if and only if $ s_D -1 \leq n \leq t_D-1$, and
 $h_B(n-s_D)$ is nonzero if and only if $s_D \leq n \leq s_D+e(B)+2$.
Since $t_D  > s_D+e(B)+2 $, we see $h_D(n)= s_D$ whenever $h_B(n-s_D)$ is nonzero ($h_B$ so to speak sits
on the top of $h_D$ as in figure \ref{fig}). Now it follows from
$h_C (n)= h_B(n-f)+h_D(n)$  that  $\lambda_C= \lambda_B \cup \lambda_D$.
\end{proof}

\bd \label{trivial}
Suppose $D_0 \subseteq C_0$ are curves in $\pso$ contained in a surface $X$, and
$D$ is obtained from $D_0$ by an elementary biliaison of height $h$ on $X$. The biliaison is defined
by an injective morphism $v: \ideal{D_0,X}(-h) \ra \coo_X$ whose image is $\ideal{D,X}$.
Then 
the image of the restriction of $v$ to $\ideal{C_0,X}(-h)$, is the ideal
$\ideal{C,X}$ of a curve $C \subset X$, obtained by biliaison from $C_0$. In this case,
we say that the biliaison from $C_0$ to $C$ is {\em induced} by the given biliaison from $D_0$ to $D$. Note that $C$ contains $D$.
\ed

\br
When $D_0$ is empty, a biliaison induced from $D_0$ is the same thing as a trivial biliaison. Indeed, in this case
$v$ is multiplication by a local equation of the complete intersection $D$ in $\coo_X$,
and $v$ maps $\ideal{C_0,X}(-h)$ onto $\ideal{C_0,X} \ideal{D,X}$.
\er

\br
For an elementary biliaison from $C_0$ to $C$ to be induced by a biliaison of $D_0$
it is enough that the corresponding morphism $u: \ideal{C_0,X}(-h) \ra \coo_X$  lift to a morphism
$\hat{u}: \ideal{D_0,X}(-h) \ra \coo_X$. Indeed, $\hat{u}$ is automatically injective because
its kernel $\mathcal{K}$ is isomorphic to a subsheaf of $\ideal{D_0,C_0}(-h) \subseteq \coo_{C_0}(-h)$, and
at the same time is a subsheaf of $\coo_X(-h)$; since  $\coo_{X}$ and $\coo_{C_0}$  have no common associated points,
we must have $\mathcal{K}=0$.
\er

\bl \label{commute}
Suppose $C_0$ contains $B$ and $D_0$, and $I_{D_0}/I_{C_0} \cong R_B(-f)$.
Suppose $C$ is obtained by an elementary biliaison from $C_0$
induced by an elementary biliaison of height $h$ from $D_0$ to $D$ on a surface $X$.
Then $C$ contains $D$ and $B$, and
$$
I_{D}/I_{C} \cong R_B (-f-h).
$$
\el
\begin{proof}
Since the biliaison from $C_0$ to $C$ is induced by that from $D_0$ to $D$, $C$ contains $D$, and
$$
I_{D}/I_{C} \cong \frac{I_{D_0}/I_X (-h)}{I_{C_0}/I_X (-h)} \cong R_B(-f-h)
$$
In particular, $R_B (-h-f)$ is an $R_C$-module, therefore $B \subseteq C$.
\end{proof}

\bl \label{e-trivial}
Suppose $C_0$ contains $B$ and $D_0$, and $I_{D_0}/I_{C_0} \cong R_B(-s_{D_0})$. If
$k$ is an integer such that
$$
k \geq \mbox{\em Max} \; (s_{D_0} + e(B) +6, e(C_0)+4),
$$
then any height one biliaison from $C_0$ to $C$ on a surface of degree $k$
is induced by a biliaison from $D_0$ to a curve $D$ such that
$$
I_{D}/I_{C} \cong R_B  (-s_D)
$$
\el
\begin{proof}
The lemma generalizes \cite[Remark 2.7.c,  p. 65]{MDP}, which treats the case $C_0=B$ and $D_0$ is empty.
The statement in this case becomes: if
$k \geq e(C_0)+6$, then every height one elementary biliaison from $C_0$ to $C$ on a surface of degree $k$
is trivial.

To prove the statement, let $X$ be the degree $k$ surface on which the biliaison from $C_0$ to $C$ is defined,
and apply $\mbox{Hom}_{\coo_X} (\; \cdot \;, \coo_X)$ to the exact sequence
$$\exact{\ideal{C_0,X}(-1)}{\ideal{D_0,X}(-1)}{\coo_B (-s_{D_0}-1)}
$$
to see $u: \ideal{C_0,X}(-1) \ra \coo_X$ lifts to $\hat{u}: \ideal{D_0,X}(-1) \ra \coo_X$
if and only if the image of $u$ in $\mbox{Ext}^1_{\coo_X} (\coo_B (-s_{D_0}-1), \coo_X) $
vanishes. Now by Serre's duality on $X$ the latter
$\mbox{Ext}$ group is dual to
$$
H^1 (X, \coo_B (k-s_{D_0}-5))
$$
which is zero because $k \geq s_{D_0} + e(B) +6$. Thus $u$ lifts to give a height one biliaison from $D_0$ to a curve $D$ inducing
the biliaison from $C_0$ to $C$. By Lemma \ref{commute} above
$I_{D}/I_{C} \cong R_B  (-s_{D_0}\!-\!1)$. Finally, since $k \geq s_{D_0} +1$, we have $s_D=s_{D_0}+1$.
\end{proof}

\bp \label{gap}
Suppose the biliaison type $\lambda_C=(k_1,  k_2, \ldots, k_s)$ of an ACM curve
$C$ has a gap at $j$.

Then $C$ contains ACM
curves $B$ and $D$ such that
$$\lambda_B=(k_1, k_2, \ldots, k_j), \;\;\; \lambda_D=(k_{j+1}, k_{j+2}, \ldots, k_s),
\;\;\; \mbox{and} \;\;\; I_{D}/I_{C} \cong R_B (-s_D).$$
Furthermore, $(B,D)$ is the unique pair of ACM curves with the above properties.
\ep
\begin{proof}
Note  that $s=s_C$. Suppose first $j=s\!-\!1$, that is,  $k_{s}\geq k_{s\!-\!1} + 3$.
Since $k_s=e(C)+3$, by Proposition \ref{LR} $C$ is obtained by a special biliaison
on a surface $X$ of degree $k_s$ from an ACM curve $B$. By definition of biliaison type,
$\lambda_B=(k_1, k_2, \ldots, k_{s\!-\!1})$. As $k_{s\!-\!1}= e(B)+3$, we see
$$k_{s}\geq k_{s\!-\!1} + 3=e(B)+6.$$
By Lemma \ref{e-trivial} the biliaison is trivial, so $C$ contains a plane section $D$ of $X$,
and $I_{D}/I_{C} \cong R_B  (-1)$. Since
$\lambda_{D}=(\deg (X)) =(k_s)$, the statement holds when $j=s\!-\!1$

We now suppose $j<s\!-\!1$ and proceed by induction on $s-j$.
By Proposition \ref{LR} $C$ is obtained by a special biliaison
on a surface $X$ of degree $k_s$ from an ACM curve $C_0$
whose biliaison type is $\lambda_0:=\lambda_{C_0}=(k_1, k_2, \ldots, k_{s\!-\!1})$.
Thus $\lambda_0$ has a gap at $j$, and $s_{C_0}=s\!-\!1$, hence
by induction $C_0$
contains ACM curves $B$ and $D_0$ such that $\lambda_B=(k_1, k_2, \ldots, k_j)$,
$\lambda_{D_0}=(k_{j+1}, k_{j+2}, \ldots, k_{s\!-\!1})$, and  $I_{D_0}/I_{C_0} \cong R_B  (-s_{D_0})$.

In particular, $s_{D_0}=s-j-1$, so that
$$k_s \geq k_{j+1}+s-j-1\geq k_j+3+s_{D_0}=e(B)+6+s_{D_0}.$$
Since $k_s=e(C)+3 \geq e(C_0)+4$,  by Lemma \ref{e-trivial} the biliaison from $C_0$ to $C$ is
induced by a biliaison from $D_0$ to a curve $D$, and $I_{D}/I_{C} \cong R_B (-s_D)$.
Finally, since $D$ is obtained from $D_0$ by a special biliaison, $D$ is ACM and
$\lambda_{D}= \lambda_{D_0} \cup (k_s)=(k_{j+1}, k_{j+2}, \ldots, k_{s})$.

It remains to prove uniqueness. Note that $s_D=s-j$ is determined by $C$, hence so is $t_D$ because
$$t_D-s_d+1 = \mbox{Min} (\lambda_D)=k_{j+1}.$$

By assumption $e(B)+3 =k_j \leq k_{j+1}-3 =t_D-s_D-2$, hence from the exact sequence
$$
\exact{\omega_D (m)}{\omega_C (m)}{\omega_B (s_D+m)}
$$
we see
$$
H^0 (\omega_D (m)) = H^0 (\omega_C (m)) \;\;\;\; \mbox{for every $m \leq 3-t_D$}.
$$

We will show that $\Omega_D=H^0_* (\omega_D)$ is generated over the polynomial ring $R=H^0_* (\PP^3)$
by its elements of degree $\leq 3-t_D$.
Taking this for granted for the moment, it follows that $\Omega_D$ is the submodule of $\Omega_C$
generated by $\ds \bigoplus_{m \leq 3-t_D } H^0 (\omega_C(m))$, hence it is determined by $C$.
But $I_D$ is the annihilator of $\Omega_D$, because $R_D$ is Cohen-Macaulay with canonical module $\Omega_D$,
hence $D$ is determined by $C$.

Since $t_D-s_D+1=k_{j+1}> 1$, the curve $D$ is contained in a unique surface $S$ of degree $s_D$,
and therefore $B$ is also determined, being the residual curve to $D=C \cap S$ in $C$.

To finish, we need to show $\Omega_D=H^0_* (\omega_D)$ is generated by its sections of degree $\leq 3-t_D$.
For this we choose a complete intersection $Y$ of type $(s_D,u)$ containing $D$ and let $E$ be the curve linked
to $D$ by $Y$. As $\Omega_D \cong I_E/I_Y (-e_Y)$ and $I_E$ is generated by its elements
of degree $\leq e(E)+3$, it is enough to show $e(Y)-t_D \geq e(E)$.

From $\omega_E(-e(Y)) \cong \ideal{D}/\ideal{Y}$ and $h^0 (\ideal{D}(t_D-1))= h^0 (\ideal{Y}(t_D-1))$,
we see $h^0 (\omega (t_D-1-e(Y)))=0$, that is
$$
t_D - e(Y) \leq -e(E)
$$
as desired.
\end{proof}

\bc \label{lines}
Let $C \subset \pso$ be an irreducible, reduced ACM curve that is
contained in a smooth surface $X$ of degree $s=s_C$. Let $t=t_C$ and
$e=e(C)$.
\begin{enumerate}
    \item[a)]
    If $h_C(e+1)=3$, $h_C(e+2)=2$, then $C$ has a unique $(e+3)$-secant line $L$,
    and every surface of degree $\leq e+2$ containing $C$ contains $L$ as well.

    \item[b)] if $h_C(t)=s-2$, $h_C(t+1)=s-3$ (so that $s \geq 3$), then $X$ contains a line $L$ that is a
    $\;(t\!-\!s\!+\!1)$-secant of $C$.
\end{enumerate}
\ec

\br
Condition $a)$ is satisfied for example when $C$ is linked to a line by a complete intersection of type $(s,t)$;
condition $b)$ is satisfied when $C$ is linked to a plane curve of degree $s\!-\!1$ by a complete intersection of type $(s,t)$.
More generally, if $\Gamma$ is linked by $(s,t)$ to $C$, then condition $a)$ says $\lambda_\Gamma$ has a gap after $1$,
so that $\Gamma$ contains a line, and condition $b)$ says $\lambda_\Gamma$ has a gap before $s\!-\!1$, so that $\Gamma$ contains a plane
curve of degree $s\!-\!1$, which is linked to a line on $X$. See proof below.

As a partial converse, we will see in the proof of Theorem \ref{gon} that,
if, for every smooth $C$ in the Hilbert scheme $A(h)$, the general surface of degree $s$ containing
$C$ contains a line, then the $h$-vector of $C$ satisfies either $a)$ or $b)$.
\er

\begin{proof}[Proof of \ref{lines}]
Since $X$ is smooth, by definition of $t$ there is surface $X_t$ of degree $t$
containing $C$ but not $X$. Thus $C$ is contained in the complete intersection $Y=X \cap X_t$.
Let $\Gamma$ the curve linked to $C$ by $Y$.
Then on $X$
$$
C \sim tH-\Gamma
$$
where $H$ denotes a plane section of $X$, and $\sim$ stands for linear equivalence.
By \cite[Corollary 5.2.19]{M}
$$
h_{\Gamma}(n)= h_Y(s+t-2-n)-h_C(s+t-2-n).
$$

\noindent {\bf Case a)} \ \  Suppose first $h(e+1)=3$ and $h(e+2)=2$. The above formula implies
$$
s_\Gamma= \mbox{Min} \{s, s+t-4-e\}.
$$
But $t \leq e+3$ because $h_C (e+3)=0$, hence $s_\Gamma= s+t-4-e$. The conditions on $h_C$ then translate
as follows:
$$
h_\Gamma (s_\Gamma)=h_\Gamma (s_\Gamma+1)=s_\Gamma -1.
$$

If $s_\Gamma=1$, this implies $\Gamma=L$ is a line.
If $s_\Gamma \geq 2$, then the condition on $h_\Gamma$
is equivalent
to $\lambda_\Gamma=(1,k_2, \ldots)$, with $k_2 \geq 4$ because $h_\Gamma (n) \geq s_\Gamma -1$ at least for
$n=s_\Gamma-2, s_\Gamma-1, s_\Gamma, s_\Gamma+1$. By Proposition \ref{gap} $\Gamma$ contains a line $L$ and an ACM curve
$D$ with $I_D/I_{\Gamma} \cong R_L (1-s_\Gamma)$.
We can treat the two cases simultaneously if we take $D$ to be the empty curve when $s_\Gamma=1$.

By Proposition \ref{gen-add}, $\Gamma=L+D$ on $X$, and $L.D=s_\Gamma-1$.
Thus
$$C.L= (tH-L-D).L=t+s-2-s_\Gamma+1=s+t-s_\Gamma-1= e+3.$$

In particular, every surface of degree $\leq e+2$ containing $C$ contains $L$ as well.
On the other hand, $C+L$ is an ACM curve, because it is linearly equivalent to $D+tH$.
Therefore $I_C/I_{C+L} \cong R_L(-C.L)=R_L(-e-3)$.
It follows that $h_{C \cup L}(n)$ and $h_C (n)$ differ
only for $n=e+3$, where their value is $1$ and $0$ respectively.
In particular, $h_{C \cup L}(e+2)= h_C (e+2)=2$ and $h_{C \cup L} (e+3)=1$,
so that by \cite[Proposition 1.5]{nollet}
the homogeneous ideal of $C \cup L$ is generated by its forms of degree $\leq e+2$, hence by the
forms in $I_C$ of degree $\leq e+2$.

Suppose now $M$ is an $(e+3)$-secant line of $C$. Then the homogeneous
ideals of $C$ and $C \cup M$ coincide in degrees $\leq e+2$. It follows that the ideal of $C \cup L$ is contained
in that of $C \cup M$, hence $C \cup L= C \cup M$ and $L=M$.
Therefore $L$ is the
unique $(e+3)$-secant of $C$.

\vspace{.3cm}
\noindent {\bf Case b)} \ \   Suppose now $h_C(t)=s-2$ and $h_C(t+1)=s-3$.
Then $h_\Gamma (s\!-\!3)=h_\Gamma (s-2)=1$ and $h_\Gamma(s\!-\!1)=0$. This implies either
$\lambda_\Gamma =(s\!-\!1)$, or $\lambda_\Gamma=(\ldots, k_{u-1}, s\!-\!1)$ with $s\!-\!1-k_{u-1} \geq 3$.
By Proposition \ref{gap} $\Gamma$ contains a plane curve $P$ of degree $s\!-\!1$ and an ACM curve
$B$ (possibly empty) such that $I_P/I_{\Gamma} \cong R_B (-1)$.

By Proposition \ref{gen-add}, $\Gamma=B+P$ on $X$, and $B.P= d_B$. Now let $L$ be the line residual to $P$
in the intersection of $X$ with the plane of $P$. Then $B.L=B.H-B.P=0$, hence
$$C.L= (tH-B-P).L= ((t-1)H-B+L).L= t-1+2-s=t\!-\!s\!+\!1 $$
\end{proof}

Given any sequence $\lambda=(k_1,  k_2, ,  \ldots, k_u)$ with $r\!-\!1$ gaps
(for any $r \geq 1$),
we can decompose
$\lambda$ uniquely as
\begin{equation}\label{decomp}
\lambda= \lambda_1 \cup \lambda_2 \cup \cdots \cup \lambda_r
\end{equation}
where each $\lambda_i$ has no gaps and, if $a_i$ and $b_i$ denote  respectively the minimum and the maximum integer
in $\lambda_{i}$, we have $a_{i+1}-b_{i} \geq 3$. We call (\ref{decomp}) the {\em gap decomposition} of $\lambda$.

\bt \label{dthm}
Let $A(\lambda)$ denote  the Hilbert
scheme parametrizing ACM curves having biliaison type $\lambda$.
If $C$ is general in $A(\lambda)$, then $C$ is reduced and for every $ f \geq e(C)+3$, there exists a smooth surface $F$ of degree
$f$ containing $C$.

If $
\lambda= \lambda_1 \cup \lambda_2 \cup \cdots \cup \lambda_r
$
 is the gap decomposition of $\lambda$, then:
\begin{enumerate}
    \item[a)] every ACM curve $C \in A(\lambda)$ contains ACM subcurves $D_i$, $i=1,2, \ldots, r$,
    such that $\lambda_{D_i}= \lambda_i$;
   \item[b)] if $C$ is general in $A(\lambda)$, then
   $$
C=D_1 \cup D_2 \cup \ldots \cup D_r
   $$
   where the $D_i$'s are distinct smooth irreducible ACM curves satisfying $\lambda_{D_i}= \lambda_i$;
   for every $1 \leq i_1 < i_2 <  \cdots < i_h \leq r$, the curve
   $$
D_{i_1} \cup D_{i_2} \cup \ldots \cup D_{i_h}
   $$
   is ACM and has biliaison type $\lambda_{i_1} \cup \lambda_{i_2} \cup \ldots \cup \lambda_{i_h}
   $.
\end{enumerate}
\et

\br
The $D_i$'s in \ref{dthm} (for $i \geq 2$) are not necessarily general in $A(\lambda_i)$: this is because they
are forced to lie on surfaces containing $D_j$ for $j< i$.
\er

\begin{proof}[Proof of Theorem \ref{dthm}]
%
Recall that by a theorem of Ellingsrud $A(\lambda)$ is irreducible (see Remark \ref{ell}).
By Proposition \ref{gap} and induction on the number of gaps
we see that for each $i$, $1 \leq i \leq r$, there are ACM curves $C_i$ and $D_i$ with the following properties
\begin{enumerate}
    \item $C_r=C$ and $C_1=D_1$;
    \item if $2 \leq i \leq r$, $C_i$ contains $C_{i-1}$ and $D_i$, and $I_{D_i}/I_{C_i}= R_{C_{i-1}} (-s_{D_i})$;
    \item $\lambda_{D_i}= \lambda_i$ for every $1 \leq i \leq r$;
    \item $\lambda_{C_i}= \lambda_{1} \cup \lambda_{2} \cup \ldots \cup \lambda_{i}$ for every $1 \leq i \leq r$.
\end{enumerate}

We claim that for every $1 \leq i_1 < i_2 <  \cdots < i_h \leq r$ there are ACM curves
$E_{i_1, i_2, \ldots, i_h} \subseteq C_{i_h}$ such that
\begin{enumerate}
    \item if $h=1$, $E_{i}=D_i$, and, if $h=r$, $E_{1,2,\ldots,r }=C$;
    \item if $2 \leq h \leq r$, $E_{i_1, i_2, \ldots, i_h}$ contains $E_{i_1, i_2, \ldots, i_{h-1}}$ and $D_{i_h}$, and
    $$I_{D_{i_h}}/E_{i_1, i_2, \ldots, i_h}= R_{E_{i_1, i_2, \ldots, i_{h-1}}} (-s_{D_{i_h}});$$
    \item $\lambda_{E_{i_1, i_2, \ldots, i_h}}= \lambda_{i_1} \cup \lambda_{i_2} \cup \ldots \cup \lambda_{i_h}$.
\end{enumerate}
We prove the statement by induction on $h$. When $h=1$ there is nothing to prove. Suppose $h>1$. By the induction hypothesis,
there is a curve $A=E_{i_1, i_2, \ldots, i_{h-1}} \subseteq C_{i_{h-1}}$ with the above properties. Let $B=C_{i_h -1}$.
By Lemma \ref{sub} below  there exists a curve $C_0 \subseteq C_{i_h}$ containing $B$ and $D_{i_h}$
such that
$I_{D_{i_h}}/I_{C_0}  \cong R_A(-s_{D_{i_h}})$. Since $A$ and $D_{i_h}$ are ACM, it follows from \ref{gen-add} that $C_0$ is ACM
as well. We define $E_{i_1, i_2, \ldots, i_h}$ to be  $C_0$. Then $E_{i_1, i_2, \ldots, i_h}$ has the required properties
(the formula for the biliaison type follows from Proposition \ref{gen-add} d)).

To see the components $D_i$ of a generic $C$ are smooth, we follow the original proof of Gruson-Peskine \cite[2.5]{GP}.
More precisely we show:
if $\lambda= \lambda_1 \cup \lambda_2 \cup \cdots \cup \lambda_r$ is the gap decomposition of $\lambda=(k_1, \ldots, k_s)$,
then there exists an ACM curve $C$ with $\lambda_C=\lambda$ satisfying the following properties:
\begin{enumerate}
    \item $C$ is contained in a smooth surface for every $f \geq k_s=e(C)+3$;
    \item $$C=D_1 \cup D_2 \cup   \cdots \cup D_r
   $$
   where the $D_i$'s are smooth irreducible ACM curves satisfying $\lambda_{D_i}= \lambda_i$; in particular, $C$ is reduced;
    \item $\omega_{D_r} (-e(D_r))$ has a section whose scheme of zeros is smooth (i.e. contains no multiple points)
\end{enumerate}

We will prove this statement by induction on $s$ as in \cite[2.5]{GP}.
For $s=1$, the statement is about plane curves and is well known (note $e(C)+3=d_C$ for a plane curve $C$).

Assume now the statement is true for $\lambda$, fix a curve $C$ with the above properties,
and consider $\lambda^+= \lambda \cup \{k_{s\!+\!1}\}$.
We have two cases to consider:

\vspace{.3cm}
\noindent {\bf Case 1)} \ \  $k_{s\!+\!1}\leq  k_s+3$. In this case $\lambda^+$ has a gap at $s$, and its gap decomposition
is
  $\lambda^+= \lambda_1 \cup \lambda_2 \cup \cdots \cup \lambda_r \cup \{k_{s\!+\!1}\}$.

By assumption $k_{s\!+\!1} \geq k_{s}+3=e(C)+6$, thus there exists a smooth surface $X$ of degree $k_{s\!+\!1}$
containing $C$. Let $D_{r+1}$ be a general plane section of $X$, and let $C^+= C \cup D_{r+1}$. Then
$D_{r+1}$ is smooth with $\lambda=(k_{s\!+\!1})$, thus $C^+$ satisfies $(2)$ with respect to $\lambda^+$.
It also satisfies $(3)$ because $\omega_{D_{r+1}} (-e (D_{r+1})) \cong \coo_{D_{r+1}}$.
By construction $C^+$ lies on the smooth surface $X$ of degree $k_{s\!+\!1}=e(C^+)+3$. The fact that
$C^+$ is  contained in a smooth surface of degree $f$, for every $f>e(C^+)+3$, follows now
from the fact that $\ideal{C^+}(e(C^+)+3)$
is generated by its global sections - see e.g. Peskine-Szpiro and \cite[Corollary 2.9]{nollet}.
Thus $C^+$ also satisfies $(1)$, and we are done in case 1.

\vspace{.3cm}
\noindent {\bf Case 2)} \ \ $k_{s\!+\!1}=k_{s}+1$ or $k_{s}+2$. In this case the gap decomposition of $\lambda^+$ is
 $$\lambda^+= \lambda_1 \cup \lambda_2 \cup \cdots \cup \lambda_{r-1}\cup \lambda^+_r $$
where $\lambda^+_r= \lambda_r \cup \{k_{s\!+\!1}\}$.

We can still find a smooth surface $X$ of degree $k_{s\!+\!1}$ containing $C$ because $k_{s\!+\!1}>e(C)+3$.
In particular $X$ contains $D_r$. The proof of \cite[2.5]{GP}
shows that the general curve $D^+_{r}$ in the linear system $D_r+H$ on $X$ is smooth with $\lambda_{D^+_{r}}=\lambda^+_r$, and
that $\omega_{D^+_{r}} (-e(D^+_{r}))$ has a section whose scheme of zeros is smooth.
Thus $C^+=D_1 \cup D_2 \cup   \cdots \cup D^+_r$ has the required properties (note $e(C^+)+3=k_{s\!+\!1}=\deg (X)$).


\end{proof}

\bl \label{sub}
Suppose  $C \subset \PP^3$ is a curve, with subcurves $B$, $D$ such that $I_{D}/I_{C}  \stackrel{\beta}{\cong} R_B (-f)$.
If $A$ is a subcurve of $B$, then there exists a unique curve $C_0$ such that
\begin{enumerate}
    \item $C_0$ is contained in $C$;
    \item $C_0$ contains $A$ and $D$, and there is an isomorphism $I_{D}/I_{C_0}  \stackrel{\alpha}{\cong} R_A(-f)$ which makes commutative the diagram
\begin{equation*} \label{diagram}\begin{array}{ccc}
    I_{D}/I_{C}  & \stackrel{\beta}{\cong} & R_B (-f) \\
             \downarrow & & \downarrow  \\
    I_{D}/I_{C_0}  & \stackrel{\alpha}{\cong}  & R_A (-f) \\
\end{array}
\end{equation*}
where the vertical arrows are induced by the inclusions $C_0 \subseteq C$ and $A \subseteq B$.
\end{enumerate}
If $A$ and $D$ have no common components, then $C_0=A \cup D$.
\el
\begin{proof}
The inclusion
$$I_{A}/I_{B} (-f) \hookrightarrow R_B(-f) \stackrel{\beta^{-1}}{\cong} I_{D}/I_{C}  \hookrightarrow R_C$$
defines an ideal $J$ in $R_C$.
Uniqueness is clear, because if such a $C_0$ exists, we must have $I_{C_0}/I_C=J$.
To show existence, let $I$ be the inverse image of $J$ in the polynomial ring $R=H^0_*(\coo_{\PP^3})$ so that
$I/I_C \cong I_{A}/I_{B} (-f)$.
The given isomorphism $I_{D}/I_{C}  \stackrel{\beta}{\cong} R_B (-f)$
induces  $I_{D}/I \stackrel{\alpha}{\cong}  R_A (-f)$, hence an exact sequence
$$
\exact{R_A(-f)}{R/I}{R_D}.
$$
From this exact sequence we see that $R/I$ has depth at least one, hence $I$ is the
saturated ideal of a subscheme $C_0 \subset C$.

By construction $I_{C_0}/I_C \cong I_{A}/I_{B} (-f)$, so that the given isomorphism $I_{D}/I_{C}  \stackrel{\beta}{\cong} R_B (-f)$
induces  $I_{D}/I_{C_0} \stackrel{\alpha}{\cong}  R_A (-f)$ with the desired properties. Finally, we
can check $C_0$ is a locally Cohen-Macaulay curve looking at the exact sequence
$$
\exact{\coo_A (-f)}{\coo_{C_0}}{\coo_D}.
$$

If $A$ and $D$ have no common components, then $C_0$ contains the union $A \cup D$. Since both $C_0$ and $A \cup D$ are locally
Cohen-Macaulay curves of degree $d_A+d_D$, they must be equal.
\end{proof}

\section{Bounds on the quadratic form $\phi(D,D)$} \label{qsection}
Let $X \subset \PP^3$ be a smooth surface of degree $s \geq 2$. Recall from section \ref{examples} the
quadratic form on $\mbox{Pic} (X)$:
$$
\phi (D,D)= (D.H)^2-(D^2)(H^2)=\det
\begin{bmatrix}
  D.H & H^2 \\
  D^2 & D.H
\end{bmatrix}
$$

In the proof of our main theorem it will be crucial to be able to bound $\phi(D,D)$ from below
in terms of the degree $d_D$ when $D$ is an ACM curve on $X$.
Note that if $D$ is a curve on $X$, then
\begin{equation} \label{phi}
\phi(D,D)=d_D^2+s(s-4)d_D-2s(g(D)-1)
\end{equation}
Thus, if we fix the degree $d_D$ and $s$, then knowing $\phi(D,D)$ is the same as knowing
the genus $g(D)$, and bounding $\phi(D,D)$ from below is the same as bounding $g(D)$ from above.
In fact, the bounds of this section can be seen as a refinement of the bounds on
the genus of an ACM curve of Gruson-Peskine \cite{GP} - see Remark \ref{remark}.
The form $\phi(D,D)$ has the advantage of being invariant if we replace $D$ with
$mH-D$ or $D+nH$, that is, it is invariant under liaison and biliaison on $X$.
Thus one can compute
$\phi(D,D)$ assuming $D$ is a minimal curve on $X$.

To compute these bounds we note that, by equation (\ref{phi}), the form  $\phi(D,D)$ for an ACM curve $D$
depends  only on the $h$-vector (or the biliaison type $\lambda$) of $D$ and on $s$. Since it is enough
to consider only minimal curves on $X$, and there only finitely many possible biliaison types $\lambda$ of
minimal curves for each $s$, our proof will proceed by a careful analysis of these $\lambda$'s.

We call a biliaison type $\lambda$ {\em $s$-minimal} if it corresponds to a minimal ACM curve on a smooth surface $X$ of degree $s$.
Since minimal is equivalent to $e+3<s$ by \ref{minimal}, the $s$-minimal $\lambda$'s are just those increasing sequences
of positive integers $\lambda= (k_1, k_2 , \ldots , k_{u})$ satisfying $k_u <s$. There are $2^{s\!-\!1}$ such possible
sequences (including the empty one), and by \ref{acm} the corresponding curves are linked by a complete intersection $(s,s)$
to curves with $s$-basic $h$-vectors. For any such $\lambda$, we let $d,g,e$ be the corresponding invariants of the associated curve $\Gamma$,
and we define
\begin{equation}
q (\lambda)= \phi(\Gamma, \Gamma)=d^2+s(s-4)d-2s(g-1).
\end{equation}
Then one verifies the formula:
\begin{equation} \label{ql}
q (\lambda) = \sum_{i=1}^{u}  k_i(s\!-\!1)(s-k_i) \; - \;
2 \sum_{1 \leq i <  j \leq u} k_i(s-k_j).
\end{equation}
See Figure \ref{table1} for the $s$-basic $h$-vectors and associated $s$-minimal biliaison types $\lambda$ for
$s=4,5$ and a few for $s=6,7,8,9$, together with the values $q$ takes on them.

\bd
Suppose $\lambda=(k_1, k_2,\ldots , k_{u})$ is $s$-minimal. Then we define
the {\em $s$-dual}  $\lambda^\prime$ of $\lambda$ to be
$$
\lambda'=(s-k_{u}, s- k_{u-1}, \cdots, s- k_1)
$$
if $\lambda \neq \emptyset$. If $\lambda= \emptyset$, then $\lambda^\prime=\emptyset$.
Note that, if $\lambda$ is the biliaison type of
an ACM curve $\Gamma$, then $\lambda'$ is the biliaison type of
a curve linked to $\Gamma$ by a complete intersection of two surfaces of degree
$s_\Gamma=u_{\lambda}$ and $s$ (cf.  section \ref{liaison}).
\ed

\bp The invariants of $\lambda'$ are:
\begin{enumerate}
    \item  $\ds u_{\lambda^\prime}=u_\lambda$
    \item  $\ds d_{\lambda^\prime}=u_\lambda s -d_\lambda$
    \item $q(\lambda^\prime)=q(\lambda)$
\end{enumerate}
\ep
\begin{proof}
The first two equalities are obvious.
The equality $q(\lambda^\prime)=q(\lambda)$ follows from (\ref{ql}), or can be
deduced from the invariance of $\phi(D,D)$ under liaison on $X$.
\end{proof}

We say that $\lambda_1=(k_1, k_2, \ldots, k_{u})$
{\em precedes} $\lambda_2=(l_1, l_2,  \ldots, l_{v})$ and write $\lambda_1<\lambda_2$
if  $k_{u} <  l_1$. In this
case, if $\lambda_2$ is $s$-minimal, then
$$\lambda_1\cup \lambda_2= (k_1, k_2,  \cdots, k_{u}, l_1,  \cdots, l_{v})$$
is also $s$-minimal. Note that $(\lambda\cup \mu)^\prime= \mu^\prime \cup \lambda^\prime$.

\bex
A plane curve of degree $k<s$ on a surface $X$ of degree $s \geq 2$ is minimal.
The corresponding $\lambda$ sequence is $\lambda=(k)$, and
$q((k))=k(s\!-\!1)(s-k)$.
\eex

\bex More generally  if $\lambda$ is the biliaison type of a complete intersection of two surfaces
    of degrees $a \leq b < s$ then
    $$
q(\lambda)= ab(s-a)(s-b)
$$
\eex
\bex \label{compute} \
 Let $\lambda=(1,  2,  \cdots, k\!-\!1, k)$ with $k<s$.
    Then $d_{\lambda}= \ds \frac{1}{2} k (k+1)$ and
    $$
q(\lambda)= d_\lambda \left(s^2-\frac{2}{3}s (2k+1)+d_{\lambda} \right)
    $$
\eex

The first statement of Proposition \ref{recursion} below determines,
once $q((k))$ is known, the function $q(\lambda)$ by induction on the number $u_{\lambda}$ of elements of $\lambda$.

\bp \label{recursion}

Suppose $\lambda <\mu$ are $s$-minimal. Then
\begin{enumerate}
    \item[a)] $q(\lambda\cup \mu)= q(\lambda)+q(\mu)-2 d_{\lambda} d_{\mu^\prime}$

    \vspace{.2cm}
    \item[b)]
If $\lambda<(k)$ and $(k+1)<\mu$, then
$$
q(\lambda \cup (k+1) \cup \mu)-q(\lambda \cup (k) \cup \mu )= (s\!-\!1)(s\!-\!1-2k)-2(d_{\mu^\prime}-d_{\lambda})
$$
\item[c)]
Suppose  $\beta$ is another $s$-minimal biliaison type, and $h$, $k$ are two integers such that
$$ \lambda <(h-1), \quad (h)<\beta <(k), \quad (k+1)<\mu.$$
Let $\delta=\lambda \cup (h) \cup \beta \cup (k) \cup \mu$
and $\epsilon=\lambda \cup (h\!-\!1) \cup \beta \cup (k\!+\!1) \cup \mu$. Then
$$
q(\delta) -
q(\epsilon)=
2s(k-h-u_\beta) \geq 2s >0
$$
\end{enumerate}
\ep

We next show that $q(\lambda)$ increases if one inserts a new integer in a sequence $\lambda$:
\bc \label{w12} \ Let $(k_1 , k_2 ,  \cdots , k_u)$ be $s$-minimal.
\begin{enumerate}
\item[a)]
If $k_{u}<k < s$, then $$q( k_1 , k_2,  \ldots, k_{u}, k) \geq q(k_1, k_2, \ldots, k_{u}) + k(s-k)^2$$
In particular,  $q(\lambda) \geq (s\!-\!1)^2$ unless $\lambda= \emptyset$.

\vspace{.2cm}
\item[b)]
If $k_i <k < k_{i+1}$, then
$$
q(k_1, k_2, \ldots, k_i, k, k_{i+1}, \ldots,  k_u) \geq q(k_1, k_2, \ldots, k_{r})+k(s-k)
$$
\end{enumerate}
\ec
\begin{proof}
Let $\lambda=(k_1, k_2, \ldots, k_{u})$
By Proposition \ref{recursion} we have
\begin{equation*}
q(\lambda \cup (k))=q(\lambda)+q(k)-2 d_{\lambda} (s-k)= q(\lambda)+ (s-k)\left(k(s\!-\!1)-2
d_{\lambda}\right).
\end{equation*}
Thus the first claim follows from
$$\ds d_{\lambda}= \sum_{1}^{r} k_i \leq \frac{1}{2} k
(k-1). $$

For the second claim,  set $\lambda=(k_1, k_2,  \ldots, k_{i})$ and $\mu=(k_{i+1}, k_{i+2}, \ldots, k_{u})$.
Using  \ref{recursion} we compute $$q(\lambda \cup ((k) \cup \mu))-q(\lambda \cup \mu)=
q((k) \cup \mu)-q(\mu)+ 2 d_\lambda( d_{\mu'}-d_{((k) \cup \mu)'})
$$
Now $d_{\mu'}-d_{((k) \cup \mu)'}=-(s-k)$ while by duality and the first claim
$$
q((k) \cup \mu)-q(\mu)=q(\mu' \cup (s-k))-q(\mu') \geq (s-k) k^2.
$$
Hence
$$q(\lambda \cup ((k) \cup \mu))-q(\lambda \cup \mu) \geq
(s-k) k^2 -2 d_{\lambda} (s-k) = (s-k)(k^2 - 2 d_{\lambda}) \geq k(s-k).
$$
where the last inequality follows from $\ds d_{\lambda}= \sum_{j=1}^{i} k_j \leq \frac{1}{2} k
(k-1). $
\end{proof}

We now prove a lower bound for $q(\lambda)$ in terms of the residue
class of $d_\lambda$ modulo $s$.

\bp \label{lbounds} \
Let $\lambda$ be $s$-minimal, of degree $d$ congruent to $f$ modulo $s$,
with $0 \leq f <  s$. Then
\begin{enumerate}
  \item[a)]
If $u_\lambda=2$ so that $\lambda=(h,k)$ with $h+k\equiv f$ ($mod \, s$),
       then
    \begin{equation*}
    q(\lambda)
   =
    \begin{cases}f(s\!-\!1)(s\!-\!f)+2h(k\!-\!1)s \mbox{\ \ \ \  \ \ \ \ \ \ \ \ \ \ \ \ \ \ if $h+k < s$}
              \\ f(s\!-\!1)(s\!-\!f)+2(s\!-\!k)(s\!-\!h\!-\!1) s \mbox{\ \ \ \ \ \ if  $h+k  \geq s$}
               \end{cases}
    \end{equation*}
    \vspace{.2cm}
   \item[b)] If $u_\lambda \geq 3$ and $s \geq 5$
    $$
    q (\lambda) \geq 2s+ m(f,s)
    $$
    where $m(f,s)$ denotes the minimum of $q(\mu)$ as $\mu$ varies among $s$-minimal biliaison types satisfying  $u_\mu=2$ and
    $ d_\mu \equiv f$ or $d_\mu \equiv s\!-\!f$ ($mod \, s$). Its value is
$$
 \begin{cases}f(s\!-\!1)(s\!-\!f)+2s(f\!-\!2)  \mbox{\ \ \ \ \ \ \ \  if $3 \leq f \leq s\!-\!f$ or if $f=s\!-\!2,s\!-\!1$}
              \\
              f(s\!-\!1)(s\!-\!f)+2s(s\!-\!f-2) \mbox{\ \ \  if  $3 \leq s\!-\!f \leq f$ or $f=0,1,2$}
               \end{cases}
$$
and is attained by $\lambda=(1,f\!-\!1)$ and $\lambda^\prime=(s\!-\!f\!+\!1,s\!-\!1)$ when
$3 \leq f \leq s\!-\!f$ or if $f=s\!-\!2,s\!-\!1$, and by
$\lambda=(1,s\!-\!f\!-\!1)$ and $\lambda^\prime=(f+1,s\!-\!1)$ when
 $3 \leq s\!-\!f \leq f$ or $f=0,1,2$.
\end{enumerate}
\ep

\begin{proof}
The first statements is a simple computation.  To prove the second statement,
note that the role of $f$ and $s\!-\!f$ is symmetric, reflecting the fact that $q(\lambda)=q(\lambda^\prime)$.
Thus we can replace $\lambda$ with $\lambda^{\prime}$ whenever convenient.
If $\lambda=(k_1 ,  k_2, \ldots, k_{r})$ and there are two indices $i<j$ such that
$k_i -1 > k_{i-1}$ and $k_{j}+1 < k_{j+1}$,
we replace $k_i$ by  $k_i -1$ and $k_j$ by $k_j+1$ to obtain a new increasing sequence
$\lambda_1$ with the same degree as $\lambda$, hence the same $f$.
Then  $ q(\lambda) \geq q(\lambda_1)+2s$ by Proposition \ref{recursion}.c. When $u_\lambda =2$, it follows that the minimum $m(f,s)$ is attained
by sequences of the form $(1,k)$ or $(h,s\!-\!1)$, as in the statement. When $u_\lambda \geq 3$,
iterating the above procedure, and passing to the dual word if necessary, we may assume $\lambda$ is one of the following sequences:

\begin{eqnarray*}
&(1, 2, \cdots, h)& \;\;\;\; 3 \leq h < s   \\
&(1, 2, \cdots, h, s\!-\!m, s\!-\!(m\!-\!1),  \cdots, s\!-\!1) & \;\;\;\; 1 \leq m \leq h, \: 2 \leq h \leq s\!-\!m\!-\!2   \\
&(1, 2, \cdots, h, k)& \;\;\;\; 2 \leq h \leq k\!-\!2   \\
&(1, 2, \cdots, h, k, s\!-\!m, s\!-\!(m\!-\!1),  \cdots, s\!-\!1)
& \;\;\;\;
m \leq h,  \: 1 \leq h \leq k\!-\!2, \: k \leq s\!-\!m\!-\!2
\end{eqnarray*}

If  $\lambda= (1, 2, \cdots, s\!-\!1)$, we replace it with $(2, \cdots, s\!-\!2)$ as
$$
q(1, 2, \cdots, s\!-\!1) > q(2, \cdots, s\!-\!2)
$$

If $h \geq 2$, we define
$$
\mu= (2, \cdots, h\!-\!1, h+1, \cdots)
$$
to be the sequence
obtained removing $1$ and $h$ from $\lambda$ and adding $h+1$.
If $h=1$, then $\lambda=(1, k, s\!-\!1)$ with $3 \leq k \leq s\!-\!3$, in which case we define $\mu=(k+1,s\!-\!1)$.

Then $d_{\mu}=d_\lambda$, $u_{\mu}=u_\lambda-1$,
hence we will be done by induction on $u_\lambda$ if we show $q(\lambda) \geq q(\mu) +2s$. By \ref{recursion}.a we can assume
$\lambda=(1, 2, \cdots, h)$ and
$\mu= (2, \cdots, h\!-\!1, h+1)$.
Then one computes $q(\lambda)-q(\mu)=2s$.
\end{proof}

\br \label{remark}
One can show that the bound $q(\lambda) \geq f (s\!-\!1)(s-f)$ is equivalent to the bound given by Gruson and Peskine (see \cite{GP})
for the genus of an ACM curve of degree $d > s(s\!-\!1)$ not lying on a surface degree $s\!-\!1$. They also show that curves of maximal genus
are linked to plane curves: in our notation this means $u_\lambda=1$ if $q(\lambda)$ attains its minimal value $f (s\!-\!1)(s-f)$.
\er

\bc \label{cbounds}
Let $\lambda$ be $s$-minimal of degree $d$ congruent to $f$ modulo $s$,
with $0 \leq f <  s$. If $u_\lambda \geq 2$, then
$$
q(\lambda) \geq
 \begin{cases}
 2s(s\!-\!2)\mbox{\ \ \ \  \ \ \ \ \ \ \ \ \ \ \ \ \ \ \ if $f=0$;} \\
3s^2-8s\!+\!1 \mbox{\ \ \ \  \ \ \ \ \ \ \ \ \ \ \ if $f=1$ or $f=s\!-\!1$;} \\
2s^2-4s+4 \mbox{\ \ \ \ \ if  $f \notin \{0,1,s\!-\!1\}$.}
\end{cases}
$$
\ec
\begin{proof}
We may assume $s \geq 5$ because the cases $s=3,4$ are easily checked directly (cf. table \ref{table1}).
If $f = 0,1$ or $s\!-\!1$, the statement follows immediately from Proposition \ref{lbounds}.
If $f \neq 0,1, s\!-\!1$, then again by Proposition \ref{lbounds}:
$$q(\lambda) \geq q(f)+2s \geq  q(2) + 2s= 2s^2-4s+4$$
\end{proof}
\bc \label{s1bound}
Suppose $s \geq 5$ and let  $\lambda$ be $s$-minimal.
Suppose $q(\lambda) \leq (s\!+\!1)^2$. Then one of the
following occurs:
\begin{enumerate}
    \item $\lambda=\emptyset$ and $q(\lambda)=0$;
    \item $\lambda=(1)$ or $\lambda=(s\!-\!1)$, and $q(\lambda)=(s\!-\!1)^2$;
    \item $5 \leq s \leq 7$ and $\lambda=(2)$ or $\lambda=(s\!-\!2)$, so that $q(\lambda)=2(s\!-\!1)(s\!-\!2)$;
    \item $s=6$ and $\lambda=(3)$, so that $q(\lambda)=3(s\!-\!1)(s\!-\!3)=45$;
    \item $s=5$ or $6$ and $\lambda=(1, s\!-\!1)$;
    \item $s=5$ and $\lambda= (1, 3)$ or $\lambda=(2, 4)$, in which case $q(\lambda)=36=(s\!+\!1)^2$;
    \item $s=5$ and $\lambda= (1, 2)$ or $\lambda=(3, 4)$, in which case $q(\lambda)=34$.
\end{enumerate}
Furthermore, if $q(\lambda) \leq (s\!-\!1)^2$, then either (1) or (2) occurs. If $(s\!-\!1)^2 <q(\lambda) \leq s^2$,
then either $s=4$ and $\lambda =(2)$ or $(1, 3)$, or $s=5$ and $\lambda =(2)$ or $(3)$.
\ec
\begin{proof}
Suppose first $\lambda=(f)$. Then $q(\lambda)=f(s\!-\!1)(s\!-\!f)$. One checks
this is bigger than $(s\!+\!1)^2$ except in the cases listed in the
statement.

Suppose now $u_\lambda \geq 2$. 
If $f=0$, then $q(\lambda) \geq 2s(s-2)$ by Corollary \ref{cbounds}, and this is bigger than $(s\!+\!1)^2$
unless $s \leq 6$. When $s=5$ or $6$, one checks by hand the only
possibility is $\lambda=(1, s\!-\!1)$.

If $f=1$ or $s\!-\!1$, the lower bound for $q(\lambda)$ is
$$
3s^2-8s\!+\!1
$$
which is bigger than $(s\!+\!1)^2$ unless $s \leq 5$. When $s=5$, one
finds the two sequences $\lambda= (1, 3)$ or $\lambda=(2, 4)$.

If $f \neq 0,1, s\!-\!1$, then $q(\lambda) \geq 2s^2-4s+4$ which
is bigger than $(s\!+\!1)^2$ unless $s \leq 5$. When $s=5$, one finds
the two sequences $\lambda= (1, 2)$ or $\lambda=(3, 4)$ for which $q(\lambda)=34$.
\end{proof}

\section{Gonality of a general ACM curve} \label{thm}
In this section we give the proof of our main result computing the gonality of a general ACM curve in $\PP^3$.

\bt \label{gon} Assume char. $\mathbb{K}=0$.
Let $C \subset \PP^3_{\mathbb{K}}$ be an irreducible, nonsingular ACM curve with $h$-vector $h$,
and let $s=s_C$, $t=t_C$, $e=e(C)$ and $g=g(C)$.  Assume
 $s \geq 4$, and $(s,d,g)$ is not one of the following:
$(4,10,11)$,  $(5,15,26)$, $(5,16,30)$, $(6,21,50)$, $(6,22,55)$, $(6,23,60)$, $(7,28,85)$, $(7,29,91)$, $(8,36,133)$.

Suppose  there is a smooth surface $X$ of degree $s$ containing $C$
with the following properties:
\begin{enumerate}
    \item the linear system $|t H -C|$ on $X$ contains a reduced curve $\Gamma$, such that the irreducible components
    $D_1, \ldots D_r$ are  ACM curves, and
    $$
\lambda_{\Gamma} = \lambda_{D_1} \cup \lambda_{D_2} \cup \cdots \cup \lambda_{D_r}
    $$
    is the gap decomposition of $\lambda_\Gamma$.
    \item The Picard group of $X$ is
$$
\mbox{Pic} (X) = \Z[H] \oplus \Z[D_1] \oplus \cdots \oplus \Z[D_r]
$$
    \item $C$ is general in its linear system on $X$.
\end{enumerate}

Then
$$\gon(C)=d-l,$$ where $l=l(C)$ is the maximum order of a multisecant of $C$. Furthermore,
with the possible exception of the values of $(s,d,g)$ listed of \ref{delta-pos}.b,
$C$ has finitely many $g^1_{d-l}$, hence its Clifford index is
$$Cliff (C)= \gon(C)-2=d-l-2.$$

More precisely:
\begin{enumerate}
    \item[a)] If $h(e+1)=3$, $h(e+2)=2$, then the gonality of $C$ is $d\!-\!e\!-\!3$ and there is unique
    pencil of minimal degree, arising from the unique $(e+3)$-secant line of $C$
   (cf. \ref{lines}).

    \item[b)] if $h(t)=s-2$, $h(t+1)=s\!-\!3$, $t>s+3$, but the condition of case $a)$ above does not occur,
    then the gonality of $C$ is $d-(t\!-\!s\!+\!1)$, and there is unique
    pencil of minimal degree, arising from the unique $(t\!-\!s\!+\!1)$-secant
    line of $C$.

  \item[c)]
    if neither cases $a)$ or $b)$ above occurs, then
the gonality of $C$ is $d\!-\!4$, and every $g^1_{d-4}$ on $C$ arises
from a $4$-secant line, unless either
\begin{enumerate}
    \item[1)] $(s,d,g)$ is in the list of \ref{delta-pos}.b, or
    \item[2)]  $s=4$, $C \in |C_0+bH|$ where $b \geq 2$ and $C_0$ has degree $4$ and arithmetic genus $1$;
 in this case $|\coo_C(b)|$ is the unique $g^1_{d-4}$ that does not arise from a $4$-secant.
\end{enumerate}
\end{enumerate}
Finally, if $C$ has a complete base point free pencil of degree $k<d-4$, then the pencil arises either from an $(e+3)$-secant line
or from a $(t\!-\!s\!+\!1)$-secant line.
\et

\br
The conditions on $h$ in $a)$ and $b)$ are not satisfied in any of the cases listed in \ref{delta-pos}.b.
\er

\begin{proof}[Proof of \ref{gon}]
The gonality of $C$ is $\leq d-4$ by Proposition \ref{Cay}
%

Suppose $\mathcal{Z}$ is a complete base point free pencil of degree $k$ on $C$, and assume  $k\leq d-4$,
unless we are in one of the cases listed in Proposition  \ref{delta-pos}.b,
for which we assume $k \leq d-5$. We will classify these pencils as follows.
By \ref{delta-pos} the bundle $\mathcal{E}$ associated to $\mathcal{Z}$ on $X$ satisfies $\Delta (\mathcal{E}) >0$,
and then by Bogomolov's Theorem \ref{bogomolov} it follows that $\mathcal{E}$ is Bogomolov unstable.
Let $\coo_X (A)$ be the line bundle that destabilizes $\mathcal{E}$.
We will show that only the following cases can occur:
\begin{enumerate}
  \item for any $h$-vector, we can have $A=-H$; then by Corollary \ref{ms} the pencil $\mathcal{Z}$ arises from a multisecant line $L$
that is not contained in $X$. Corollary \ref{no5s}  shows that $k=\deg \mathcal{Z} =d-4$ and that there
is a finite set of such pencils.
  \item when $h(e+1)=3$ and $h(e+2)=2$, then $C$ has a unique $(e+3)$-secant line $L$, and $\mathcal{Z}=\mathcal{Z}(L)$.
  In this case $L \subset X$ and $A=L-H$.
    \item if $t>s+3$, $h(t)=s\!-\!2$, $h(t+1)=s\!-\!3$, then $C$ has a unique $(t\!-\!s\!+\!1)$-secant line $L$, and $\mathcal{Z}=\mathcal{Z}(L)$.
  In this case $L \subset X$ and $A=L-H$.
  \item $s=4$, $C \in |C_0+bH|$ where $b \geq 2$ and $C_0$ has degree $4$ and arithmetic genus $1$. In this case
 $\mathcal{Z}=|\coo_C (b)|$ and   $A=-C_0$.
  In particular, $\deg \mathcal{Z}= d-4$ and $\mathcal{Z}$ does not arise from a multisecant.
\end{enumerate}
The statement of the theorem clearly follows from this classification. For the Clifford index, we use the
fact proved by Coppens and Martens \cite{cm} that $Cliff (C)= \gon(C)-2$
when $C$ has a finite number of pencils of minimal degree.

We now proceed to classify the possible base point free complete pencils $\mathcal{Z}$ of degree $\leq d-4$.
Let $A$ be the divisor that destabilizes the bundle $\mathcal{E}$ associated to $\mathcal{Z}$. Recall that
$A$ sits in an exact sequence
$$
\exact{\coo_X(A)}{\mathcal{E}}{\ideal{W,X}(B)}
$$
where $W$ is zero dimensional and $(A-B).H >0$. From the exact sequence we see
$A-B=2A+C$ and
$$
(2A+C)^2 =(A-B)^2 \geq \Delta(\mathcal{E})=C^2 -4 k.
$$
By Proposition \ref{a2} we also have $(-A).H >0$ and $A^2 \geq 0$.

To be able to work effectively with the above inequalities, we write
$x=A.H$ for the degree of $A$, and consider as in section \ref{examples} the bilinear form on $Pic(X)$
$$
\phi(D,E)= (D.H) \, (E.H) - s \, (D.E)=
\det
\begin{bmatrix}
  D.H & H^2 \\
  D.E & E.H
\end{bmatrix}.
$$

We then obtain the following numerical restraints on $x$:
\begin{eqnarray}
\label{axx}   && -\frac{d}{2} < x < 0          \\
\label{ax2x}   && x^2\geq \phi(A,A) \\
\label{axdx}   && x^2+dx +ks \geq \phi(A,A+C)
\end{eqnarray}
the last two inequalities being equivalent
to $A^2 \geq 0$ and $(2A+C)^2 \geq C^2 -4 k$ respectively.

In $Pic (X)$ we can write $A= \sum a_i D_i+cH$ with $a_i \in \Z$, $c \in \Z$.
We wish to show $$\phi(A,A+C)\geq 0.$$
We first prove $\phi(D_i,D_j)<0$.
Let
$\lambda_{\Gamma}= \lambda_1 \cup \lambda_2 \cup \cdots \cup \lambda_r$
be the gap decomposition of $\lambda_\Gamma$, so that $\lambda_{D_i}=\lambda_i$.
If $i<j$,
$D_i+D_j$ is ACM with $\lambda_{D_i+D_j}=\lambda_i \cup \lambda_j$ by Theorem \ref{dthm}.
Since $\phi(D,D)= q(\lambda_D)$ for an ACM curve $D$ with $s_D<s$, by  Proposition \ref{recursion}
\begin{equation}\label{dij}
\phi(D_i,D_j)= -d_{\lambda_i}d_{\lambda'_j}<0
\end{equation}

(note that the formula $\phi(D_i,D_j)= -d_{\lambda_i}d_{\lambda'_j}$ is correct only for $i<j$).

To simplify notation we let $q_i=\phi(D_i,D_i)$ and
 $\ds b_i = - \sum_{j \neq i} \phi(D_i,D_j)$. We claim that for every $i$
$$
q_i > 2b_i
$$
To prove this let $\ds E_i= \sum_{j \neq i} D_j$. Then
$$
\phi(\Gamma,\Gamma)=\phi(D_i+E_i, D_i+E_i)= \phi(D_i,D_i)+ \phi(E_i, E_i)+2\phi(D_i, E_i)
=\phi(E_i, E_i) +q_i -2 b_i
$$
thus it is enough to show $\phi(\Gamma,\Gamma)>\phi(E_i, E_i)$, that is
$$
q(\lambda_{\Gamma}) > q(\lambda_{E_i}).
$$
The latter inequality holds by Corollary \ref{w12}, hence $q_i > 2b_i$.

We now compute
\begin{eqnarray}
\nonumber \phi(A,A)&=&
\sum_i a_i^2 \phi(D_i,D_i)+ 2 \sum_{ i < j} a_i a_j \phi(D_i,D_j)=
\\
\nonumber
&=&
\sum_i a_i^2 (q_i-b_i)- \sum_i a_i^2 \sum_{j \neq i} \phi(D_i,D_j)
+ 2 \sum_{i < j} a_i a_j \phi(D_i,D_j)=
\\
\nonumber
&=& \sum_i a_i^2 (q_i-b_i)- \sum_{ i < j} (a_i-a_j)^2 \phi(D_i,D_j) =\\
\nonumber \phi(A,C)&=& \phi (\sum_i a_i D_i, t_C H-\sum_j D_j )= \phi (\sum_i a_i D_i, -\sum_j D_j )
\\
\nonumber
&=&  - \sum_{i,j} a_i \phi(D_i,D_j) = - \sum_i a_i (q_i-b_i)
\end{eqnarray}

Therefore
\begin{eqnarray}
\label{qA}
\phi(A,A)
&=& \sum_i a_i^2 (q_i-b_i)- \sum_{ i < j} (a_i-a_j)^2 \phi(D_i,D_j) \\
\label{AC}
\phi(A,C)&=&  - \sum_i a_i (q_i-b_i)
\\
\label{QAC}
\phi(A,A+C) &=&
\sum_i (a_i^2-a_i) (q_i-b_i)- \sum_{ i < j} (a_i-a_j)^2 \phi(D_i,D_j)
\end{eqnarray}

The last equality implies $\phi(A,A+C)\geq 0$ because the $a_i$ are integers, $q_i>2b_i \geq b_i$ and
$\phi(D_i,D_j)<0$.

We now show that $\phi(A,A+C)\geq 0$ implies $x \geq -s\!-\!1$.

By hypothesis $k \leq d-4$, therefore
 $$ x^2+dx+(d-4)s \geq x^2+dx + k s \geq \phi(A,A+C) \geq 0.$$
Let $\delta$ be the discriminant of the equation
$x^2+dx +(d-4) s = 0$:
$$
\delta= d^2-4sd+16s= (d-2s)^2-4s(s-4).
$$
Let $y=d-2s$.
Since $C$ is ACM and $s=s_C$, we have  $d \geq \frac{1}{2} s(s\!+\!1)$ by \ref{degreebound}, hence
 $$y-2= d-2s-2 \geq \frac{1}{2} (s^2-3s-4) \geq  \frac{1}{2} (s^2-4s).$$
In fact, we can have equality only if $s=4$ and $d=10$, while
the hypotheses of the Theorem when $s=4$ require $d$ to be at least $11$.
Thus $y-2> \frac{1}{2} s (s-4)$ and
$$
\delta= y^2-4s(s-4) > y^2-8y+16 =(y-4)^2
$$
Thus $\delta$ is positive, and the equation has two real roots,
one smaller than $\ds-\frac{d}{2}$, the other one, say $\bar{x}$, larger than $\ds-\frac{d}{2}$.
Since $ \ds  -\frac{d}{2}< x < 0$, we conclude $x \geq \bar{x}$. Furthermore,
unless $s=4$ and $d=11$, under the hypotheses of the theorem $y-4 \geq 0$, hence
\begin{eqnarray*}
 \bar{x} &=& -\frac{d}{2}+ \frac{1}{2} \sqrt{\delta} \\
   &> & -\frac{d}{2}+ \frac{1}{2} \sqrt{y^2-8y+16} \\
   &=& -\frac{d}{2}+ \frac{1}{2} (y-4)=-s-2
\end{eqnarray*}
The inequality $ \bar{x} > -6$  holds also in case $s=4$ and $d=11$.
Thus $x \geq -s\!-\!1$.
Then from $x^2\geq \phi(A,A)$ we see
$$(s\!+\!1)^2 \geq \phi(A,A).$$

If all the $a_i$'s are zero, then $A=cH$ (note that this is the case if $C$ is a complete intersection
of $X$ and another surface). Since $-s\!-\!1 \leq x= \deg A <0$, we must have $A=-H$.

If not all the $a_i$'s are zero, let
$1 \leq i_1 < \cdots <i_h \leq r$ be the indices for which $a_i \neq 0$.
Formula (\ref{qA}) holds with this new set of indices, and shows that, if all the coefficient $a_i$'s are non zero,
then $\phi(A,A)$ attains its minimum
when all the $a_i$'s are equal to $1$. Thus
$$
\phi(A,A) \geq \phi(D,D)
$$
where $D=D_{i_1}+ \cdots +D_{i_h}$ is the support of $A$.

Now $D$ is ACM with biliaison type  $\lambda_D= \lambda_{i_1}  \cup \cdots \cup  \lambda_{i_h}$
by Theorem \ref{dthm}. If $\lambda_D$ is not
one of the special cases listed in Corollary \ref{s1bound}, then
$$
\phi(D,D)=q(\lambda_D) >(s\!+\!1)^2$$
contradicting
$(s\!+\!1)^2 \geq \phi(A,A)$.

Suppose now $\lambda_D $ is
one of the special cases listed in Corollary \ref{s1bound}. We still have
$\phi(A,A) \geq (s\!-\!1)^2$ because $\lambda_D$ is not empty.
Before examining the various cases, let us remark that, if only one of the $a_i$'s is nonzero,
so that $$A=aD+cH$$
with $D$ irreducible and $a \neq 0$, then
either $a=1$ or $a=-1$. This follows from
$$
a^2=\frac{\phi(A,A)}{\phi(D,D)} \leq \frac{(s\!+\!1)^2}{(s\!-\!1)^2}<4
$$
Also note that $D$ is irreducible precisely when $\lambda_D$ has no gaps, that is,
in all cases of \ref{s1bound} except when $s=5$ or $6$ and $\lambda= (1, s\!-\!1)$.

\begin{figure}[htb]
\caption{$s$-basic $h$-vectors and $s$-minimal biliaison types} \label{table1}
\begin{center}
\includegraphics[scale=0.8]{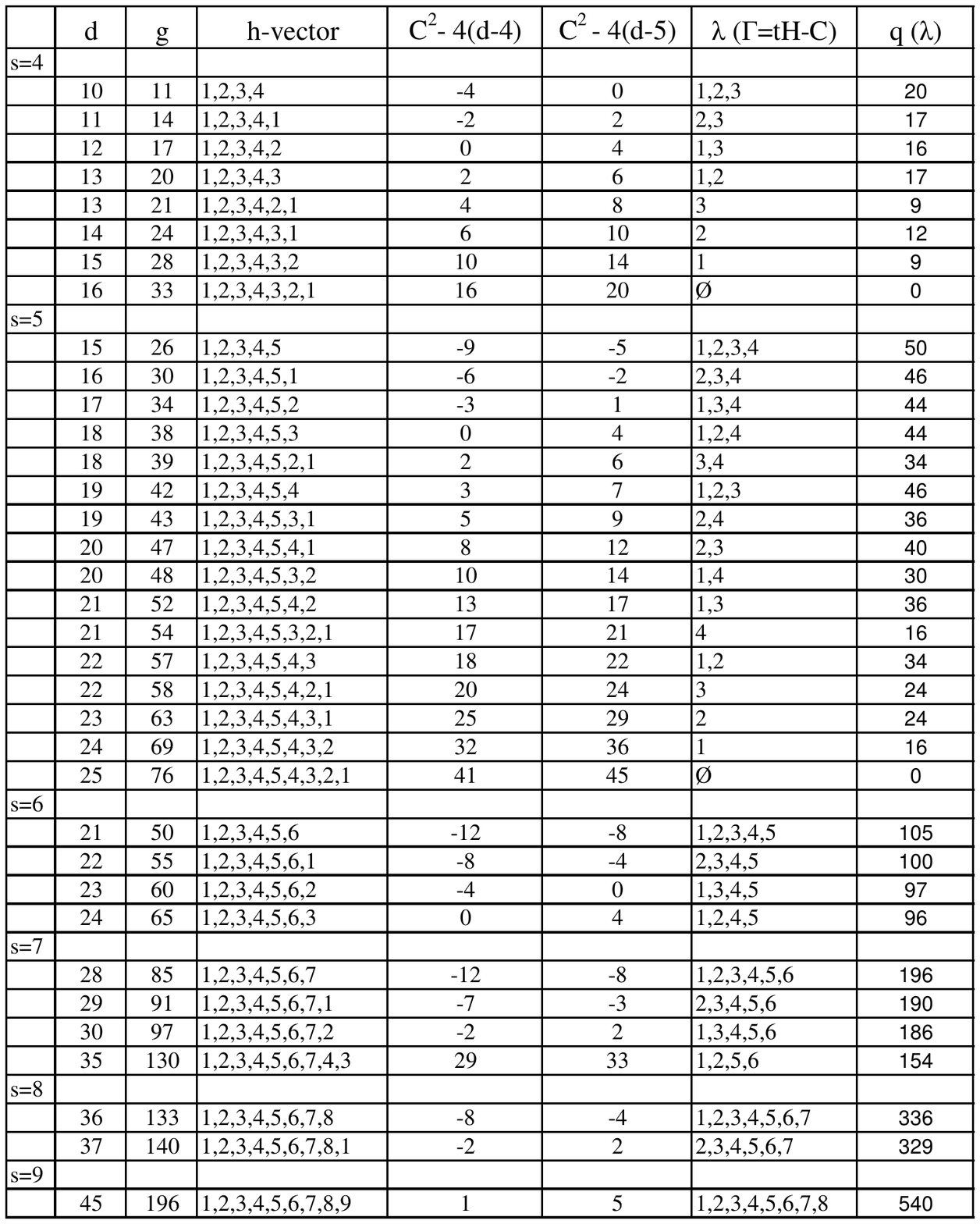}
\end{center}
\end{figure}

To complete the list of \ref{s1bound}, observe that for $s=4$ (cf. Table \ref{table1}) there are $7$ possibilities for
$\lambda_D$, because $\lambda \neq \emptyset$ and $u_{\lambda}<4$, namely
$$
(1),(2),(3),(1,2),(1,3),(2,3),(1,2,3).
$$

\vspace{.2cm}
\noindent \textbf{Case 1: } \ assume $\lambda_D \neq (1)$, $\lambda_D \neq (s\!-\!1)$, and,
when $s=5$ or $6$,  $\lambda_D \neq  (1, s\!-\!1)$.

\vspace{.2cm}
Then  $\phi(D,D) > (s\!-\!1)^2$ and $\lambda_D$ has no gaps by \ref{s1bound}. Thus $D$ is irreducible, $A=aD+cH$
with $a = \pm 1$ and
$$(s\!+\!1)^2 \geq x^2 \geq \phi(A,A)=a^2 \phi(D,D) > (s\!-\!1)^2.$$
Hence  $x=-s\!-\!1$ or $ x =-s$.


\vspace{.2cm}
\noindent \textbf{Case 1.a: \ $a=1$, $x=-s\!-\!1$} \

\vspace{.2cm}
In this case $d_D \equiv x \equiv -1$ ($mod \, s$), and by \ref{s1bound}
we must have  $ s \leq 5$. Furthermore by (\ref{axdx})
$$x^2+d x +(d-4)s \geq 0,$$ that is
$$
s^2+2s\!+\!1-sd-d + (d-4)s \geq 0
$$
so $d \leq s^2-2s\!+\!1$. This gives $d \leq 9$
if $s=4$, and $d\leq 16$ if $s=5$, while $d \geq \frac{1}{2}s(s\!+\!1)$
because $C$ is an ACM curve $s_C=s$. Thus we must have $s=5$, and
examining the list in \ref{s1bound} we find $\lambda_D=(1,3)$ is the only possibility.
Then, for $\Gamma= t H-C$, we know $\lambda_\Gamma$ contains $\lambda_{D}=(1,3)$ in its
gap decomposition and $u_{\lambda_{\Gamma}}<5$. This forces $\lambda_\Gamma=\lambda_D$,
hence $D=\Gamma$ and therefore
$$d= st -\deg (\Gamma) \geq 25-4=21$$
a contradiction, so this case does not occur.

\vspace{.2cm}
\noindent \textbf{Case 1.b: \ $a=1$, $x=-s$} \

\vspace{.2cm}
In this case $d_{D} \equiv x \equiv 0$ ($mod \, s$) and $s^2=x^2 \geq q(\lambda)$. By
\ref{s1bound} the only possibility is $s=4$ and $\lambda_D=(1, 3)$, which forces $D=\Gamma=tH-C$.
Furthermore, we must have $\gon(C)=k=d-4$ for (\ref{axdx}) to hold.

Since $x=-4= \deg (D+cH)$, we see $c=-2$. Now pick an effective divisor $C_0 \in |-A|=|2H-D|$. Then
$C_0$ is ACM with biliaison type $(1,3)$, thus $C_0$
is up to a deformation with constant cohomology an elliptic quartic.
By construction $C \in |C_0+bH|$ with $b=t-2 \geq 2$
(note that $b=2$ gives $(d,g)=(12,17)$ which is in the list \ref{delta-pos}.b).
For $b \geq 2$ the restriction of $|C_0|$ to $C$ is $|\coo_C(b)|$, and is a  $g^1_{d-4}$ on $C$
that does not arise from a multisecant.

\vspace{.2cm}
\noindent \textbf{Case 1.c:} \ $a=-1$, $x=-s\!-\!1$  or $-s$

\vspace{.2cm}
In this case $A=-D+cH$, hence, if $D=D_i$,
$$\phi(A,A)+ \phi(A,C) =
2\phi(D_i,D_i)+\sum_{j \neq i} \phi( -D_i, -D_j) = 2q_i-b_i \geq \frac{3}{2} q_i > \frac{3}{2} (s\!-\!1)^2
$$
Therefore
$$x^2+d x+(d-4) s \geq \frac{3}{2} (s\!-\!1)^2$$
which contradicts both $x =-s\!-\!1$ and $x=-s$, so this case does not occur.

\vspace{.2cm}
\noindent \textbf{Case 2: } \ $\lambda_D =(1)$, so that $D$ is a line $L \subset X$, and $A=cH+aL$ with $a = \pm 1$.

\vspace{.2cm}
In this case either $\Gamma=L$ and $\lambda_\Gamma=(1)$, or $\lambda_\Gamma$ has a gap at the beginning:
$$
\lambda_\Gamma= (1,4, \ldots)
$$
In both cases $L=D_1$  is unique.
The proof of Corollary \ref{lines} shows that the $h$-vector of $C$ satisfies $h_C(e+1)=3$ and $h_C(e+2)=2$,
and that  $C.L=e+3$. Thus in any case
$$\deg (Z)= \gon(C) \leq d-e-3.$$

We wish to show that $A=L-H$ and $Z=\mathcal{Z}(L)$.

Recall that the degree $x$ of $A$ must satisfy the inequalities $-s\!-\!1 < x <0 $ and
$$
x^2 \geq a^2 \phi(L,L)=(s\!-\!1)^2.$$
We also know $x=cs+a$ with $a=\pm 1$.
Therefore $c=-1$ and either $A=-H-L$ or $A=-H+L$.

Suppose first $A=-H-L$. Since $\deg(X)= s\geq 4$,
$$
H^0 \coo_X (H+L) \cong H^0 \coo_X (H)
$$
thus every curve $B$ in the linear system $|-A|=|H+L|$ contains the line $L$. This
contradicts Proposition \ref{a2}, according to which we can find two effective divisors in
$|-A|$ meeting properly. So $A=-H-L$ is impossible.

Therefore $A=-H+L$, and $Z=\mathcal{Z}(L)$  by Corollary \ref{ms}.
%

\vspace{.2cm}
\noindent \textbf{Case 3: } \  $\lambda_D =(s\!-\!1)$, so that $D=H-L$ is a plane curve of degree $s\!-\!1$,
residual to a line $L$ in a plane section of $X$. Furthermore, $A=cH+aD=(c+a)H-aL$ with $a=\pm 1$.\

\vspace{.2cm}

In this case $D=D_r$, thus $L$ is unique, and
either $\Gamma=D_r$ or $\lambda_\Gamma$ has a gap at the end.  The proof of Corollary \ref{lines} shows that the $h$-vector of $C$ satisfies
$h_C(t)=s-2$, $h_C(t+1)=s\!-\!3$ and that  $L$  is a $(t\!-\!s\!+\!1)$-secant line for $C$.
An argument analogous to the one of the previous case
shows $A=-H+L$, so that $\mathcal{Z}=\mathcal{Z}(L)$.

\noindent
\textbf{Case 4: } \ $\lambda_D=(1, s\!-\!1)$ with $s=5$ or $6$, hence
$A=cH+a_1 L_1 +a_2 P$ where $L_1$ is a line, $P$ is a plane curve of degree $s\!-\!1$,
and $a_1$ and $a_2$ are non zero.
Note that $\phi( L_1, P)  = - 1$, therefore
$$
\phi(A,A)= (a_1^2+a_2^2) (s\!-\!1)^2 - 2 a_1 a_2
= (a_1^2+a_2^2) (s^2-2s) +(a_1-a_2)^2 \geq 2(s^2-2s) > s^2
$$
On the other hand, $(s\!+\!1)^2 \geq x^2 \geq \phi(A,A)$. Therefore we must have
$x=-s\!-\!1$ and $a_1^2+a_2^2 <3$, that is, $a_1$ and $a_2$ can only be $1$ or $-1$.

Then
$$
-s\!-\!1=x = cs+ a_1+a_2(s\!-\!1)
$$
from which we see $-1 \equiv a_1-a_2$ ($mod \, s$). This is impossible because $a_1= \pm1$ and $a_2=\pm1$.

This complete the list of possible cases, and proves the classification of complete base point free pencils
$\mathcal{Z}$ of degree $ \leq d-4$, hence the theorem

\end{proof}

\br
In the first of the cases excluded in the theorem, namely
$s=4$ and $(d,g)=(10,11)$, we can prove $\gon (C)=6=d-4$ by the method of \cite{mjm}.
\er

\bt \label{general}
Assume the ground field is the complex numbers. Then the conclusions of theorem \ref{gon} hold for
the general ACM curve $C$ in $A(h)$.
\et

\begin{proof}
Since the conclusions of theorem \ref{gon}  are semicontinuous on $A(h)$ (cf. \cite{AC}), it is enough to show
the existence of a single curve $C$ for which the hypothesis of \ref{gon} are satisfied.

To check this,
let $h'$ denote the $h$-vector of a curve $\Gamma$ linked  by two surfaces of degrees
$s$ and $t=t_C$ to $C \in A(h)$.

Note that $h'$ may not be of decreasing type, but in any case $s_{\Gamma} \leq e_{\Gamma}+3 < s$ by
Lemma \ref{decreasing}. By  Theorem \ref{dthm}
a general curve $\Gamma$ in $A(h')$ is reduced, its irreducible components are ACM, with biliaison
type prescribed by $\lambda_\Gamma$; and, since $s > e_{\Gamma}+3$, there exist smooth surfaces of any degree
$\geq s-1$ containing $\Gamma$.

Now let $h_2$ be the $h$-vector of a curve $C_2$ linked to $\Gamma$ by the complete intersection of
two smooth surfaces of degree $s\!-\!1$ and $s$ respectively.
The flag Hilbert schemes parametrizing pairs $(\Gamma,Y)$, where $\Gamma \in A(h')$ and $Y$ is a complete intersection
of type $(s\!-\!1,s)$, is irreducible \cite[VII\S 3]{MDP}.
Thus a general $\Gamma$ in $A(h')$ can be linked to a general $C_2 \in A(h_2)$.
By Lemma \ref{decreasing} $h_2$ is of decreasing type, hence we may assume $C_2$ is smooth, and lies on smooth
surfaces of degree $s\!-\!1$ and $s$.
Since we are working over the complex numbers, we can use the Noether-Lefschetz type theorem of
Lopez \cite[Theorem II 3.1]{lopez}. We  apply this theorem to $C_2$ with $d=s$, $e=1$, and
$T$ a smooth surface of degree $s\!-\!1$ through $C_2$ to conclude that,
if $X$ is a {\em very general} surface of degree $s$ containing $C_2$, then $Pic(X)$ is freely generated by the
classes of a plane section $H$ and of the irreducible components of $\Gamma$
(here very general means, as usual, outside a countable union of proper subvarieties).

Now on $X$ we can take for $C$ a general curve in the linear system $$|C_2+(t-s+1)H|=|tH-\Gamma|.$$
The hypotheses of \ref{gon} are then satisfied for the smooth surface $X$ and the curve $C$.

One can simplify the argument using a more recent result of Brevik and Nollet \cite[Theorem 1.1]{bs}
that allows one to work directly with $\Gamma$ rather than $C_2$.
\end{proof}

\newpage

\bigskip


\end{document}